\documentstyle[aps,amstex,pre,multirow,preprint,lscape,epsfig,subfigure,hhline]{revtex}

\def\f{\frac}
\def\l{\left}
\def\r{\right}
\def\p{\partial}

\begin{document}

\title{\Large \bf A practical spectral method for hyperbolic conservation laws}

\author{Yu-Hui Sun$^{1}$, Y. C. Zhou$^{1}$  and  G. W. Wei$^{1,2}$\footnote{Corresponding author,
wei@@math.msu.edu}}
\address{$^1$Department of Mathematics,
     Michigan State University,
     East Lansing, MI 48824\\
$^2$ Department of Computational Science, National University of
Singapore, Singapore 117543}

\maketitle


\date{\today}

\begin{abstract}
A class of high-order lowpass filters, the discrete singular convolution (DSC)
filters, is utilized to facilitate the Fourier pseudospectral method for
the solution of hyperbolic conservation law systems. The DSC filters are
implemented directly in the Fourier domain (i.e., windowed Fourier pseudospectral 
method), while a physical domain algorithm is also given to enable the 
treatment of some special boundary conditions.
By adjusting the effective wavenumber region of the DSC filter,
the Gibbs oscillations can be removed effectively while the high resolution
feature of the spectral method can be retained. The utility and effectiveness
of the present approach is validated by extensive numerical experiments.

\end{abstract}


\newpage

\section{Introduction}

It is well known that due to their accuracy and efficiency, spectral
methods have great advantages over local methods in solving applicable
scientific and engineering problems \cite{Boyd,CHQZ,Fornberg,GHO,Trefethen}.
Given a hyperbolic system of nonlinear conservation laws
\begin{eqnarray}\label{CL}
\textbf{u}_t + \textbf{f}(\textbf{u})_x = 0,
\end{eqnarray}
with initial condition
\begin{eqnarray}
\textbf{u}(x,0)=\textbf{u}_0(x),
\end{eqnarray}
its solution may not exist in the classical sense because of
possible discontinuities. The direct use of spectral methods to
this problem will encounter the Gibbs oscillations
\cite{gottlieb97}, which lead to blow-up in the time integration.
Therefore, it has been of great interest in making spectral
methods applicable to the hyperbolic conservation law systems in
the past two decades. The expectation, as point out by  Gottlieb
\cite{gottlieb94}, is that the use of spectral methods will enable
us not merely to capture the shock, but also to capture the
delicate features and structures of the flow. Obviously, spectral
approaches will be extremely valuable to a class of aerodynamical
problems that involve the interaction of both turbulence and shock
\cite{LLM}. Such an interaction is some of the most challenging
problems in computation fluid dynamics.

The main objective of previous studies is to recover smooth
solutions from those that are contaminated by Gibbs oscillations
and, meanwhile, to improve the rate of convergence. For
shock-capturing, many up-to-date local methods, such as weighted
essentially non-oscillatory (WENO) scheme \cite{LOC,qiu,Shu1},
central schemes \cite{Bianco,Kurganov,LiuTad,Tadmor},
arbitrary-order non-oscillatory advection scheme \cite{toro98},
Gas kinetic \cite{Xu01}, and image processing based schemes
\cite{Sonar,Weicpc02}, perform well. The success of these local
shock-capturing schemes lies in their appropriate amount of
intrinsic numerical dissipation, which is introduced either by
explicit artificial viscosity terms, or by upwinding, or by
appropriate local average (non-oscillatory central schemes), or by
relaxation \cite{jin96}. The characteristic decomposition based on
Roe's mean matrix can also be considered as a local averaging on
the Jacobian matrix $\textrm{A(u)}
=\p\textrm{{f(u)}}/\p\textrm{{u}}$. In his recent work, LeVeque
 \cite{leveque} outlined the relation between some approximate Riemann
solvers and relaxation schemes. The relative amount of numerical dissipations
may explain why the local characteristic decomposition is not necessary in
low-order methods while it seems to be indispensable in high-order
methods \cite{qiu}.
However, when these local methods are used in the cases where flow
structures with fine details are needed to be resolved together with
shocks, numerical dissipations are usually found to be so large that they smear the
fine details \cite{LLM}. The spectral methods, on the contrary, contribute very little
numerical dissipation and dispersion in principle when used to approximate spatial
derivatives. Nevertheless, when they are applied to the approximation of spatial
derivatives on a domain containing discontinuities, the Gibbs oscillations
again have to be suppressed by appropriate means, for which there are two general
routines: (1) explicit artificial viscosity, e.g. spectral viscosity (SV)
method proposed by Tadmor \cite{tadmor89} and (2) filtering, which is the
central issue of the present study.

At present, most filters are constructed in the spectral domain and they can
be called spectral filters. Hussaini {\it et al.} \cite{hussaini84} suggested
a list of  typical spectral filers, i.e. Lanzos filter,
raised cosine filter, sharpened raised cosine filer and
exponential cutoff filter. More sophisticate and effective filters
in spectral domain are Vandeven's $p$-th order filter
 \cite{vandeven91}, Cai {\it et al.}'s sawtooth function
 \cite{cai89} and Gottlieb and Tadmor's regularized Dirichlet
function \cite{gottliebtadmor84}.

Apart from methods that are implemented in the spectral domain,
filters in the physical domain are also developed. To design an
appropriate filter in the physical domain is generally more
difficult than in the Fourier domain. A straightforward procedure
is to make use of numerical dissipation partially contained in {\it
high-order} shock-capturing schemes \cite{cai93}, such as the ENO scheme,
where, actually, numerical dissipations were introduced both in
the Fourier domain (via an exponential filter) and in the physical
domain (via ENO polynomial interpolation). Such a strategy was
generalized by Yee {\it et al.} \cite{Yee} to the so-called characteristic
filter method (with finite difference or finite volume) which has
been successfully applied to shock-capturing and turbulence
simulations. Another procedure was initiated by Gottlieb
 \cite{gottlieb92} by using the Gegenbauer polynomial to
resolve the oscillatory partial Fourier summation. Some promising
numerical results from the filter approach can be found in Refs.
\cite{wai94,Don98,gottlieb01,hesthaven99}.

Nevertheless, some prominent issues may hinder the success of
filtered spectral methods for practical numerical computations involving
shocks. In our view, there are two very important issues that are vita to the
success of a filtered spectral scheme. The first issue is how to
select optimal filters for shock-capturing.
The issue is very complex and cannot have a unique answer at this point
because there are so many properties to mind, such as flatness, ripple, filter
length, effective frequency range and length of transition band, to name only a few.
Loosely speaking, it is desirable to have filters that are free of dispersion
errors, flat while having very small transition band, short in length while
having high resolution. Moreover, what adds to the complexity is that solutions
to different hyperbolic conservation law systems may have different Fourier spectral
distributions. Consequently, one faces the difficulty that an optimal filter for
a problem may not even work  for another one. Therefore, it is desirable to have
a filter which accounts for this change by an adjustable parameter.

The second issue also concerns the change of solution, but is due to the
time evolution for a given hyperbolic conservation law system. As the Fourier
spectral distribution of a given problem changes with time, the question
is how to implement the filter? It will be too dissipative for many systems
if a filter is applied at each time step of the integration. Therefore, an
adaptive implementation, which is  controlled by a sensor during the time
integration, is appropriate. Given the complexity, it is unlikely that
there will be an ideal solution to these issues in the immediate future.
Obviously, these challenging and interesting problems call for the further
study of spectral filter approaches.

The objective of this work is to construct, analyze and implement
a class of high-order lowpass filters which can be used either in
the Fourier domain or in the physical domain. These filters are
constructed within the framework of the discrete singular
convolution (DSC) algorithm \cite{WPWjcp02,Weijcp99,WZXijnme02}
and they are used to suppress  Gibbs oscillations arising in  the
Fourier pseudospectral method (FPM). Such lowpass filters have a
wide effective wavenumber range that makes it possible to capture
the fine flow structures, which is a desirable objective of
spectral methods for hyperbolic conservation law systems. It is
this effective wavenumber range that controls the resolution of
the overall method. In our design, this effective range can be
varied according to the resolution requirement by a parameter.
Different DSC kernels have different magnitude responses and
adjustability in the Fourier domain, which in turn influence the
accuracy and resolution. The DSC filters used in this work were
originally designed in a conjugate filter oscillation reduction
(CFOR) \cite{guwei,weigu,ZhoWeijcp} scheme. By `conjugate filters'
it means that the effective wavenumber range of the lowpass filter
is largely overlapped with that of the highpass  filter which is
used for the approximation of spatial derivatives. The resulting
Fourier domain algorithm is essentially a windowed  Fourier 
pseudospectral method for shock-capturing.

The rest of this paper is organized as follows. In Section II, a brief introduction
will be devoted to the construction of the DSC lowpass filters followed by a
discussion of two types of implementations of these filters to spectral methods.
Extensive numerical experiments are carried out in Section III, including
Burgers' equation, shock-tube problems, shock-entropy wave interaction,
shock-vortex interaction. The shock-capturing ability and the high resolution feature
of the present method are addressed. Concluding remarks end this paper.

\section{Theory and algorithm}

\subsection{DSC lowpass filter}

In the context of distribution theory, a singular convolution is defined by
\begin{eqnarray}
f(x) = (T * \eta)(x)=\int_{-\infty}^{\infty} T(x-t)\eta(t) dt,
\label{eq:sc}
\end{eqnarray}
where $T(x)$ is a singular kernel and $\eta (x)$ is an element of
the space of test functions. Interesting singular kernels include
that of Hilbert type, Abel type and delta type. The former two
play important roles in the theory of analytical functions,
processing of analytical signals, theory of linear responses and
Radon transform. Since delta type kernels are the key element in
the theory of approximation and the numerical solution of
differential equations, we focus on  singular kernels of the
delta type
\begin{eqnarray}
T(x) = \delta^{(q)}(x) \; , ~~ q=0,1,2,\cdots \; , \label{eq:delta_k}
\end{eqnarray}
where superscript $(q)$ denotes the $q$th-order `derivative' of the delta
distribution $\delta (x)$, with respect to $x$, which should be understood as
generalized derivatives of distributions. When $q=0$, the kernel $T(x)=\delta(x)$
is important for the interpolation of functions. In this work,
only the case of $q = 0$ will be involved. One has
to find appropriate approximations to the above singular kernel, which can not be
directly realized in computers. To this end, we consider a
sequence of approximations
\begin{eqnarray}
\lim_{\alpha \rightarrow \alpha_0}
\delta_{\alpha}^{(q)}(x) = \delta^{(q)}(x) \; , ~~ q=0,1,2,\cdots \label{eq:seq}
\end{eqnarray}
where $\alpha$ is a parameter which characterizes the approximation
with the $\alpha_0$ being a generalized limit. Among various candidates of
approximation kernels \cite{Weijcp99}, regularized Shannon kernel (RSK)
\begin{eqnarray}
\delta_{\sigma,\Delta}(x-x_k)=\frac{ \sin \frac{\pi}{\Delta}(x-x_k)}
{\frac{\pi}{\Delta}(x-x_k)}
\exp \left( -\frac{(x-x_k)^2}{2 \sigma^2} \right) \label{eq:rsk}
\end{eqnarray}
is frequently used due to its close link to the Shannon sampling theorem.
In this formula, $\Delta$ is the grid spacing and $\sigma$ determines the width of the
Gaussian envelop. For a given $\sigma \neq 0$, the limit of $\Delta \rightarrow 0$
reproduces the delta kernel (distribution). With the RSK,
a function $u$ can be approximated by a discrete convolution
\begin{eqnarray}
u(x) \approx \sum_{k=-W}^{W} \delta_{\sigma,\Delta}(x-x_k) u(x_k) \;,
~~ q=0,1,2, \cdots \;, \label{eq:dsc}
\end{eqnarray}
where $\{x_k\}_{k=-W}^{W}$ are a set of discrete grid points which are centered
at $x$, and $2W+1$ is the computational bandwidth, or effective kernel support,
which is usually smaller than the computational bandwidth of the spectral
method - the entire domain span. Equation  (\ref{eq:dsc}) is often referred as
a DSC and the RSK is used as the prototype of a DSC lowpass filter.

\begin{figure}
\centering
\includegraphics[width=8cm]{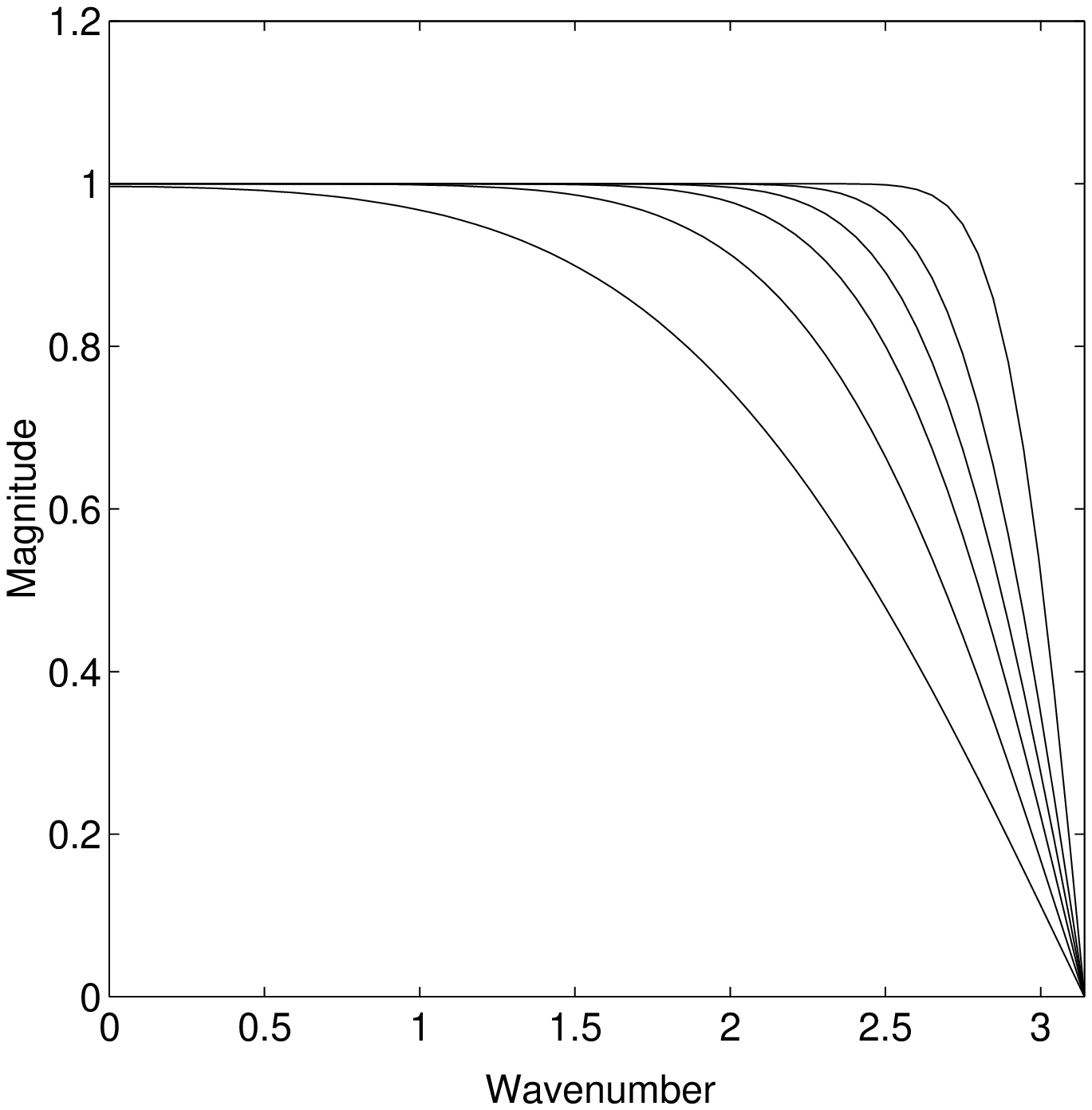}
\vspace{2mm}
\caption{Magnitude response of the DSC-RSK filter.
$r=1.0, 1.5, 2.0, 2.5, 3.2, 5.0$ from the left to right.}
\label{dsckernel}
\end{figure}

The feasibility of the DSC-RSK as lowpass filters in numerical computation is based on
an important observation: by adjusting the ratio $r$ ($r=\f{\sigma}{\Delta}$) of the
regularizer, one can control the effective wavenumber region and thus,  control the
resolution of the lowpass filter. We plot the frequency response profile of the DSC-RSK
lowpass filter in FIG. \ref{dsckernel}. A distinct property of the RSK lowpass
filter is that there is no ripples in its passband; thus the application of
this filter will not impart the low frequency. In contrast, some lowpass filters
constructed with certain polynomials may have ripples in their passband, which
would rule out their application to numerical integrations, especially for
long time evolution, where the solution may be distorted after repeated filtering.

The maximum effective wavenumber (frequency) range of the DSC-RSK is determined by
the ratio $r$, which in turn is determined by computational bandwidth $W$. The
philosophy is simple since with the increase of $W$, the Gaussian window should
also be enlarged to fully utilize the stencil. The optimal values $r$ at different
$W$ can be estimated and they characterize the effective wavenumber range. The larger
the optimal $r$ is, the wider the effective wavenumber range will be.  When
$r$ value is chosen larger than the optimal value for a given $W$, the accuracy
of the resulting filter will be degraded by oscillations. Hence we do not
recommend the use of $r$ beyond the optimal value unless in some circumstances
where the filer is utilized to stabilize a long-time integration.

\subsection{Fourier domain algorithm}

We denote $\zeta(x_j)=\delta_{\sigma,\Delta}(x_j)$ the proposed DSC-RSK lowpass filter in the
physical domain $[0,L]$ and $\hat {\zeta}(\omega_n)$  its image in the 
Fourier domain, where $\Delta=L/N$, $x_j=j\Delta$ and $\omega_n={2\pi n\over L}$.
The obvious relation between $\zeta(x_j)$ and $\hat {\zeta}(\omega_n)$ is
\begin{eqnarray}
\hat{\zeta}(\omega_n) = \Delta \sum_{j=0}^{N-1} {\zeta}(x_j) e^{-i\omega_n x_j}.
\end{eqnarray}
As in the previous work \cite{gottlieb94}, filtering in the Fourier
domain is usually accomplished by multiplying the given Fourier coefficients
$\hat{u}(\omega_n)=\Delta \sum_{j=0}^{N-1} {u}(x_j) e^{-i\omega_n x_j} $ 
by $\hat {\zeta} (\omega_n)$ directly, i.e., a windowed Fourier transform,
\begin{eqnarray}
   u^{\zeta} (x_j) = {1\over L}\sum_{n=-{N\over2}}^{{N\over2}-1} \hat{\zeta}(\omega_n) 
    \hat{u}(\omega_n) e^{i\omega_n x_j}.
\label{eq:filter1}
\end{eqnarray}
Therefore, the {\it modified} Fourier coefficients $\hat{\zeta} (\omega_n)
\hat{u}(\omega_n)$ in Eq. (\ref{eq:filter1}) are the Fourier coefficients
of a {\it localized} version of the original function $u(x_j)$: $u^{\zeta} (x_j)$.
Hence, errors arising from the discontinuity are also localized and the
accuracy away from the discontinuity can be ensured. In practical applications,
one only needs to carry out the filtering and differentiation $D_xu=u_x$ in one step
\begin{eqnarray}
   u_x^{\zeta} (x_j) = {1\over L}\sum_{n=-{N\over2}}^{{N\over2}-1} i\omega_n\hat{\zeta}(\omega_n) 
   \hat{u}(\omega_n) e^{i\omega_n x_j}.
\label{eq:filter1-2}
\end{eqnarray}
A standard
forth-order Runge-Kutta (RK-4) scheme is employed for the time advance scheme in this
work. We use a TVD switch to adaptively activate  the application of the DSC-RSK
lowpass filter.  Whenever the total variation of the approximate solution in the
two consecutive time steps exceeds a prescribed criteria, the lowpass filter is
activated. The DSC-RSK lowpass filter is only applied at the end of each
 full RK-4 cycle. The reader is referred to our earlier work 
\cite{guwei,weigu,ZhoWeijcp} for the detail of this adaptive filtering.

\subsection{Physical domain algorithm}

It is more difficult to apply the filter in the physical domain than in the
Fourier domain. On one hand, it is not easy to quantitatively analyze the
influence of numerical dissipation in the physical domain. On the other hand,
the convolution operation
\begin{eqnarray}
u^{\zeta} (x) = \int u(x-\xi) \zeta(\xi) d \xi,
\end{eqnarray}
complicates the physical domain application of the filter. However,
in some practical applications, it turns out to be quite necessary to
implement the filter in the physical domain, especially for problems
whose boundaries require some special treatments. The reflective boundary
is one of such examples that need to be treated by using grid points
outside the computational domain. In this work, we propose the method
to apply our DSC-RSK lowpass filter in the physical domain. Due to the
fact that the RSK lowpass filter is an interpolatory kernel, we consider
a two-step procedure. Step one is the prediction of mid-mesh values
from the nodal values
\begin{eqnarray}
u^\zeta_{i+\f{1}{2}}=\sum^{-1}_{j=-W} u_{i+j+1} \zeta((j+\f{1}{2})\Delta)
              +\sum^{W}_{j=1} u_{i+j}\zeta((j-\f{1}{2})\Delta).
\end{eqnarray}
The other is the reconstruction of the nodal values from the mid-mesh values
\begin{eqnarray}
  u^{\zeta\zeta}_{i}=\sum^{-1}_{j=-W} u^\zeta_{i+j+\f{1}{2}}\zeta((j+\f{1}{2})\Delta)
        +\sum^{W}_{j=1} u^\zeta_{i+j-\f{1}{2}} \zeta((j-\f{1}{2})\Delta),
\end{eqnarray}
where, $u_i$ denotes $u(x_i)$ and $\Delta$ is the grid spacing. The lowpass
filter with a $r$ value smaller than the optimal value is used in one of
these two steps. The activation of the  RSK lowpass filter
is controlled in the same manner as that of the frequency domain algorithm.

As described, the DSC-RSK lowpass filter can be applied either in the
physical domain or in the Fourier domain.  In general, filtering in the
Fourier domain is much easier and more convenient than in the physical
domain because no additional summation is required to implement
the filter in the Fourier spectral method. In contrast, the physical
domain implementation requires two additional summations.
Therefore, we have used the filter in the
Fourier domain in most of our numerical experiments except for the cases
in which the reflective boundary condition is involved. We have checked
that whenever applicable, results obtained from the Fourier domain approach
and the physical domain approach are identical.

\subsection{Post processing filter}

\begin{figure}
\centering
\includegraphics[width=8cm]{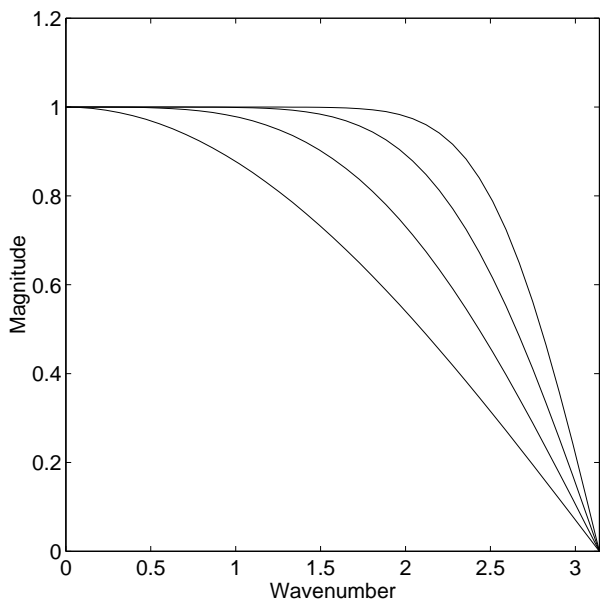}
\vspace{2mm}
\caption{Magnitude response of the Lagrange filter.
Lagrange-2, Lagrange-4, Lagrange-8 and Lagrange-16 from the left to right.}
\label{lagkernel}
\end{figure}

Once the time integration has been completed, a strong filter is usually
required for post processing the solution so as to make the solution more
presentable. In principle, we only need to reduce the $r$ value of the
DSC-RSK lowpass filter at the last time step of the computation.
As an alternative, we provide a series of lowpass
filters based on the Lagrange delta kernel which is given by
\begin{eqnarray}
\delta_{W}(x,x_j)=\prod_{k=-W, k\ne j}^{W} \frac {x-x_k}{x_j-x_k}, ~~~~ j=-W,\cdots, W.
\label{eq:lagker}
\end{eqnarray}
FIG. \ref{lagkernel} shows the magnitude response of the Lagrange
delta kernel. In fact, the highpass versions of these filters are graphically
identical to centered finite difference (FD) schemes. In our numerical experiments,
we adopt Lagrange-2 or Lagrange-4 as post processing filters. Their
efficiency is demonstrated by our numerical results.

\section{Numerical experiments}

In this section, we examine the validity and demonstrate the performance of
proposed windowed (filtered) Fourier spectral scheme for conservation law
systems. A large collection of standard
linear and nonlinear benchmark problems \cite{Shu1}, including
the linear advection equation, 1D inviscid Burgers' equation,
shock tube (Sod and Lax) problems, 1D shock-entropy interaction,
2D shock-entropy interaction, 2D shock-vortex interaction,
2D advection of an isentropic vortex and flow past a cylinder, are considered
in the present study.
Some of these problems have periodic boundary conditions therefore the
windowed Fourier pseudospectral method will be applied directly.
For problems whose boundary conditions are not periodic,
we symmetrically double the computational domain in both
$x$ and $y$ directions to compute the flux derivative.
It is noted that the extension is necessary for constructing periodic
problem feasible for the Fourier spectral method.
Our filtering procedures, as described earlier,
can be applied with general boundary conditions, such as Dirichlet or Neumann.
In the appendix, we list the filter parameter $r$ used for each numerical example.

\subsection{Scalar conservation law systems}

We begin our numerical experiments with the scalar conservation
law system, which is given by
\begin{eqnarray}
u_{t}+f(u)_{x}=0,
\end{eqnarray}
where $f(u)$ is a function of $u$. It is generally believed that
spectral methods are not very suitable for simple shock profiles
and should not be used if one is interested primarily in the
shock profiles rather than other detailed structures of the flow.
In this study, we choose scalar conservation law systems to illustrate
that the proposed global spectral scheme performs well for a class of
localized shock problems.

\textit{\underline {Example 1.}} We first solve a linear advection
equation, which is of the form
\begin{eqnarray}
   \begin{array}{rl}
      u_{t}+u_{x}= &  0 \quad  -1<x<1, \\
      u(x, 0)=  & u_0(x) \quad  \text{periodic,}
   \end{array}
\end{eqnarray}

where, $u_0$ is an initial value
\begin{eqnarray}
u_0(x)=
\l \{
   \begin{array}{ll}
      \f{1}{6} (G(x, \beta, z-\delta)+G(x, \beta, z+\delta)+4G(x, \beta, z))
      &  \quad  -0.8 \le x \le -0.6 \\
      1                  &   \quad  -0.4 \le x \le -0.2 \\
      1-|10(x-0.1)|      &   \quad   0   \le x \le  0.2 \\
      \f{1}{6} (F(x, \alpha, a-\delta)+F(x, \alpha, a+\delta)+4F(x, \alpha, a))
      &  \quad  0.4 \le x \le 0.6 \\
      0                  &    \quad   \text{otherwise}.
   \end{array}
\r.
\end{eqnarray}
Here, functions $G$ and $F$ are defined as
\begin{eqnarray}
  \begin{array}{rl}
   G(x, \beta, z)= & e^{-\beta (x-z)^2}  \\
   F(x, \alpha, a)= & \sqrt{\mbox{max}(1-\alpha^2 (x-a)^2,0)},
  \end{array}
\end{eqnarray}
where $a=0.5$, $z=-0.7$, $\delta=0.005$, $\alpha=10$ and
$\beta=\f{\text{log}2}{36 \delta^2}.$

\begin{figure}
\centering
    \includegraphics[width=9cm]{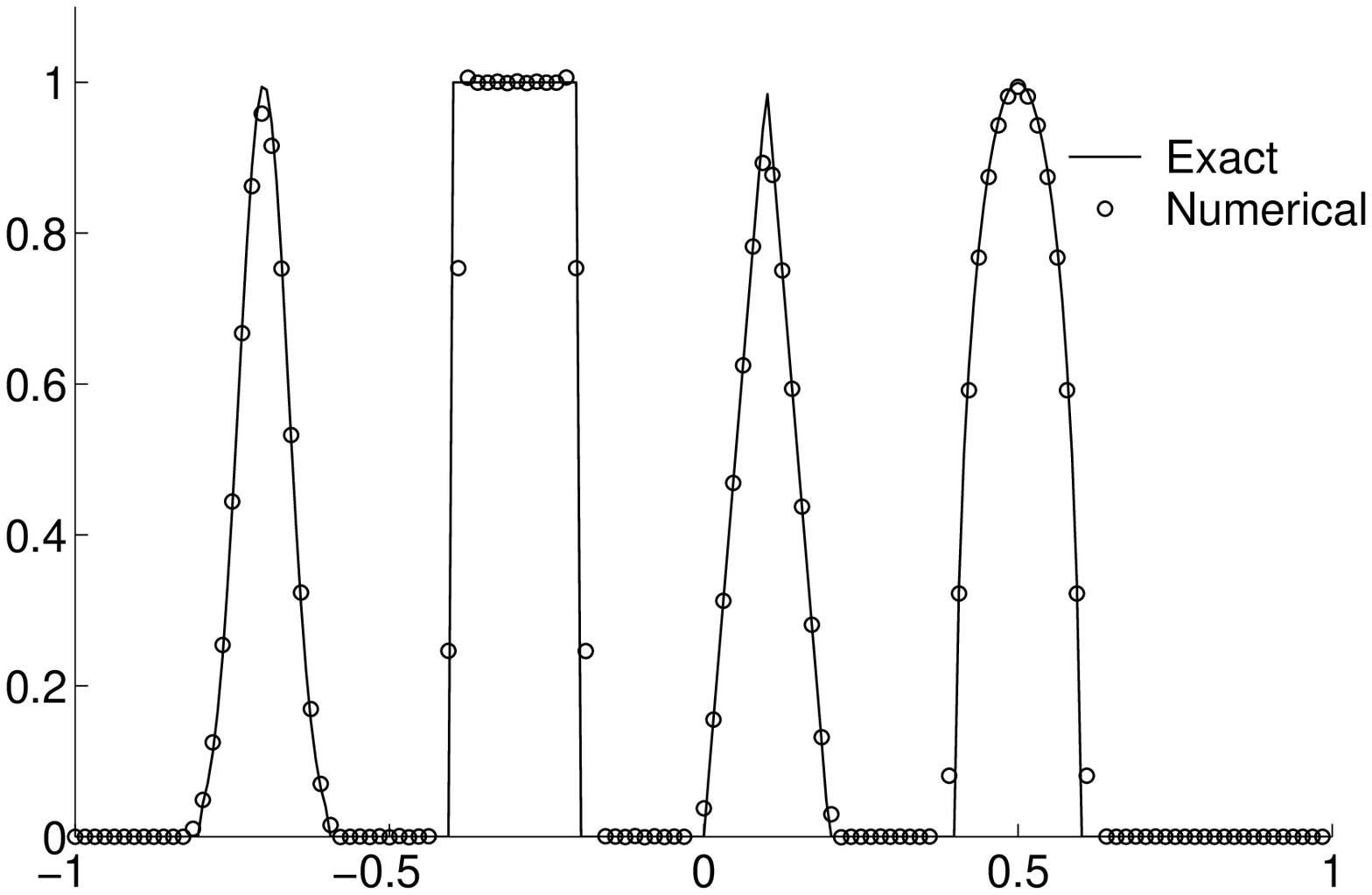}
    \caption{Linear advection equation. $t$=8, $\Delta$t=0.001, 128 grid points.}
    \label{linear1}
\end{figure}

\begin{figure}
\centering
    \includegraphics[width=9cm]{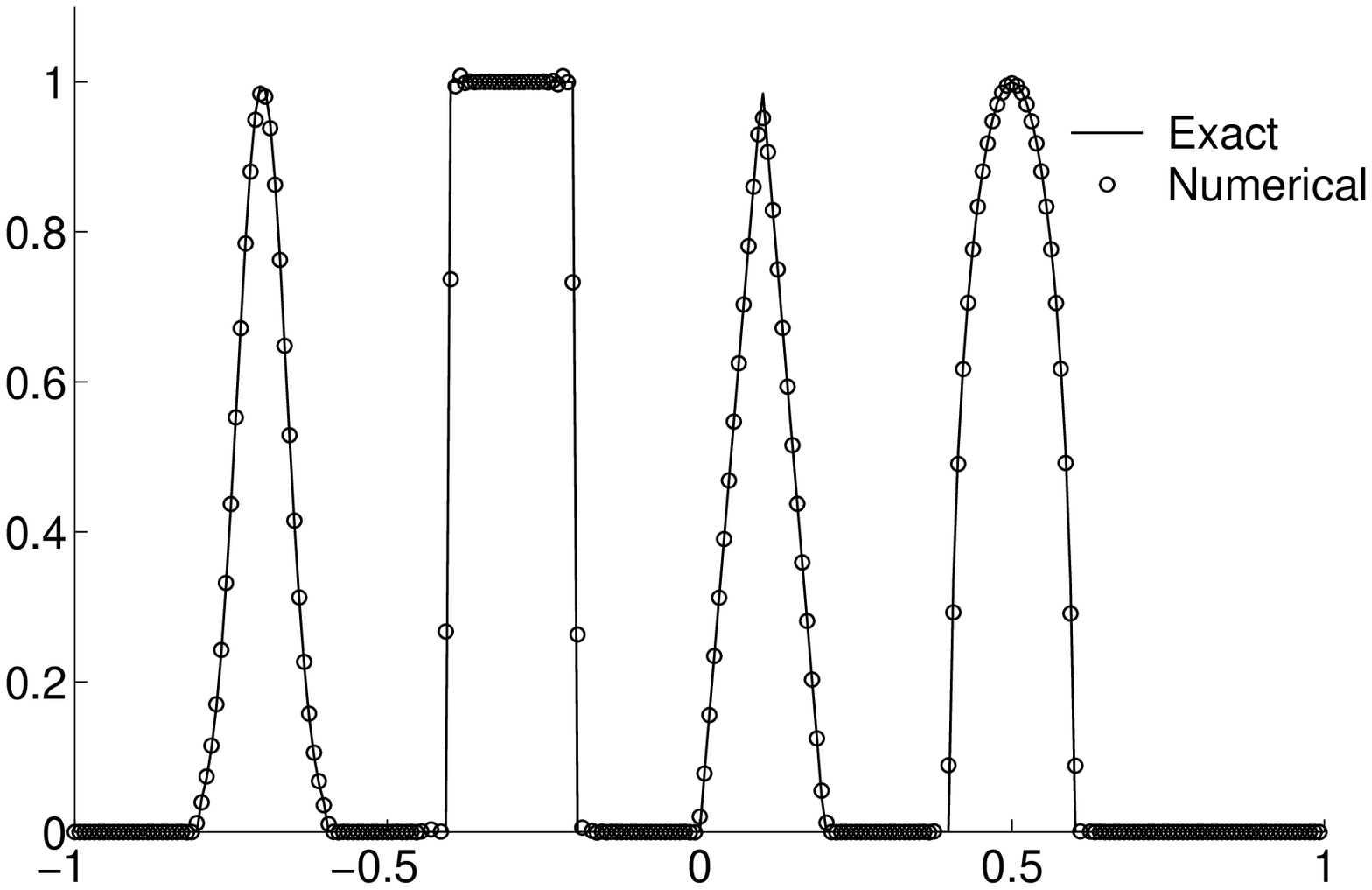}
    \caption{Linear advection equation. $t$=8, $\Delta$t=0.001, 256 grid points.}
    \label{linear2}
\end{figure}

The solution of this problem contains a combination of a Gaussian,
a square wave, a sharp triangle wave and a half ellipse. The exact
solution at any time can be easily obtained by a translation of the
initial solution at  speed 1. In FIG. \ref{linear1} and FIG. \ref{linear2},
we display the numerical results for this problem at $t$=8 with
128 and 256 grid points, respectively. We observe that our method
works well in both cases.

\begin{figure}
\centering
    \includegraphics[width=9cm]{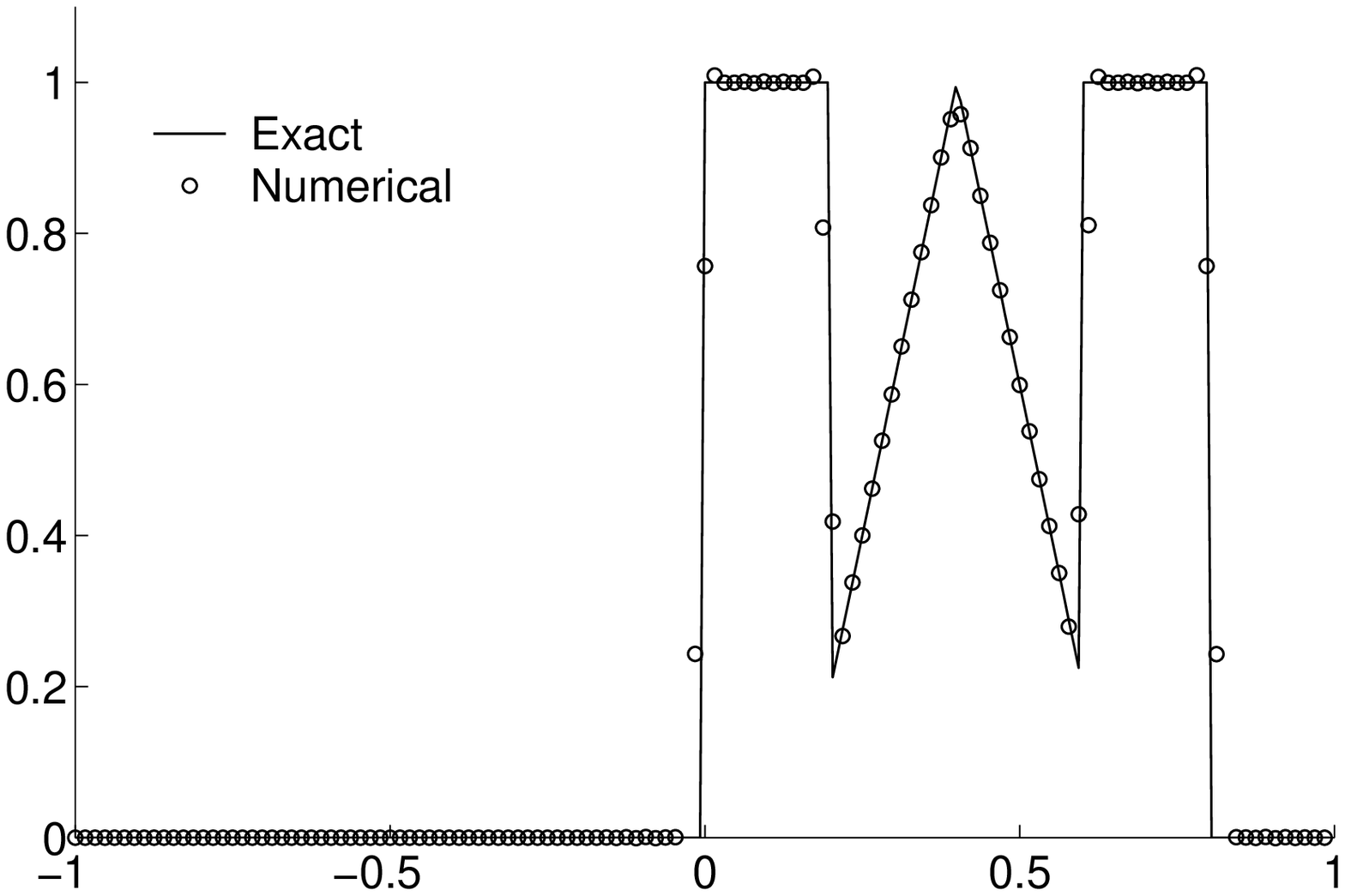}
    \caption{Linear advection equation. $t$=8, $\Delta$t=0.001, 128 grid points.}
    \label{linear3}
\end{figure}

\begin{figure}
\centering
    \includegraphics[width=9cm]{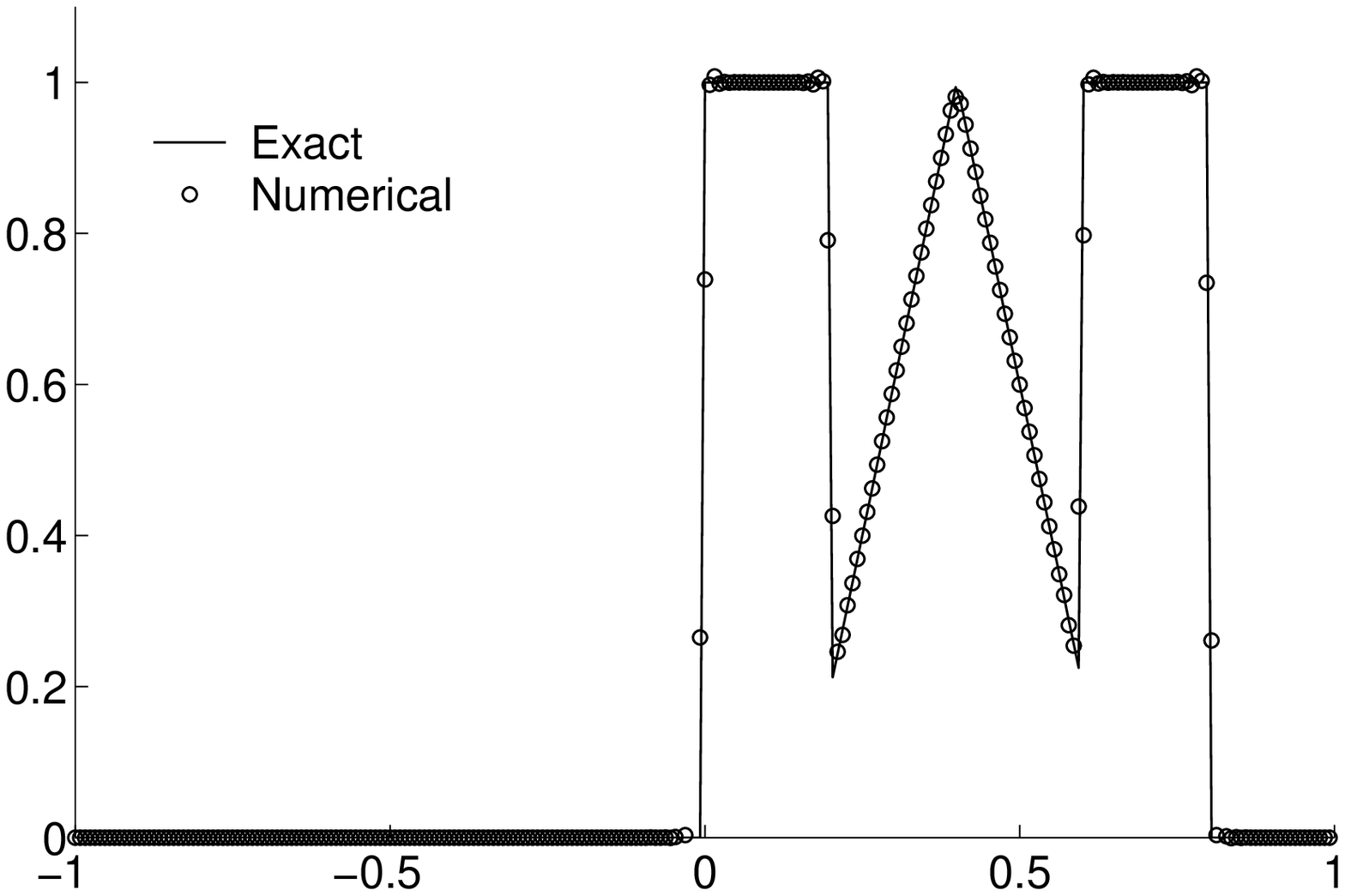}
    \caption{Linear advection equation. $t$=8, $\Delta$t=0.001, 256 grid points.}
    \label{linear4}
\end{figure}

\textit{\underline {Example 2.}}  Secondly, we test the proposed method by
considering a moving W-shape wave, i.e., a piecewise continuous initial value,
for the linear advection equation with an initial value \cite{weigu}
\begin{eqnarray}
u_0(x)=
\l \{
   \begin{array}{ll}
      1               &   \quad  0    \le x \le  0.2 \\
      4x-\f{3}{5}    &   \quad  0.2  \le x \le  0.4 \\
      -4x+\f{13}{5}  &   \quad  0.4  \le x \le  0.6 \\
      1               &   \quad  0.6  \le x \le  0.8 \\
      0               &   \quad  \text{otherwise}.
   \end{array}
\r.
\end{eqnarray}
This case is very similar to Example 1. It contains the so-called
contact discontinuity and is quite difficult to solve in hyperbolic
conservation laws. In particular, the interaction of two lines at $x=0.2$
and $x=0.6$ is difficult to resolve. Hence, it is a good test for
the shock-capturing ability of the present method. We compute
the solution up to $t$=8. The results
are shown in FIG. \ref{linear3} and FIG. \ref{linear4}.

\begin{figure}
 \centering
 \subfigure[ ]{
     \includegraphics[width=7cm]{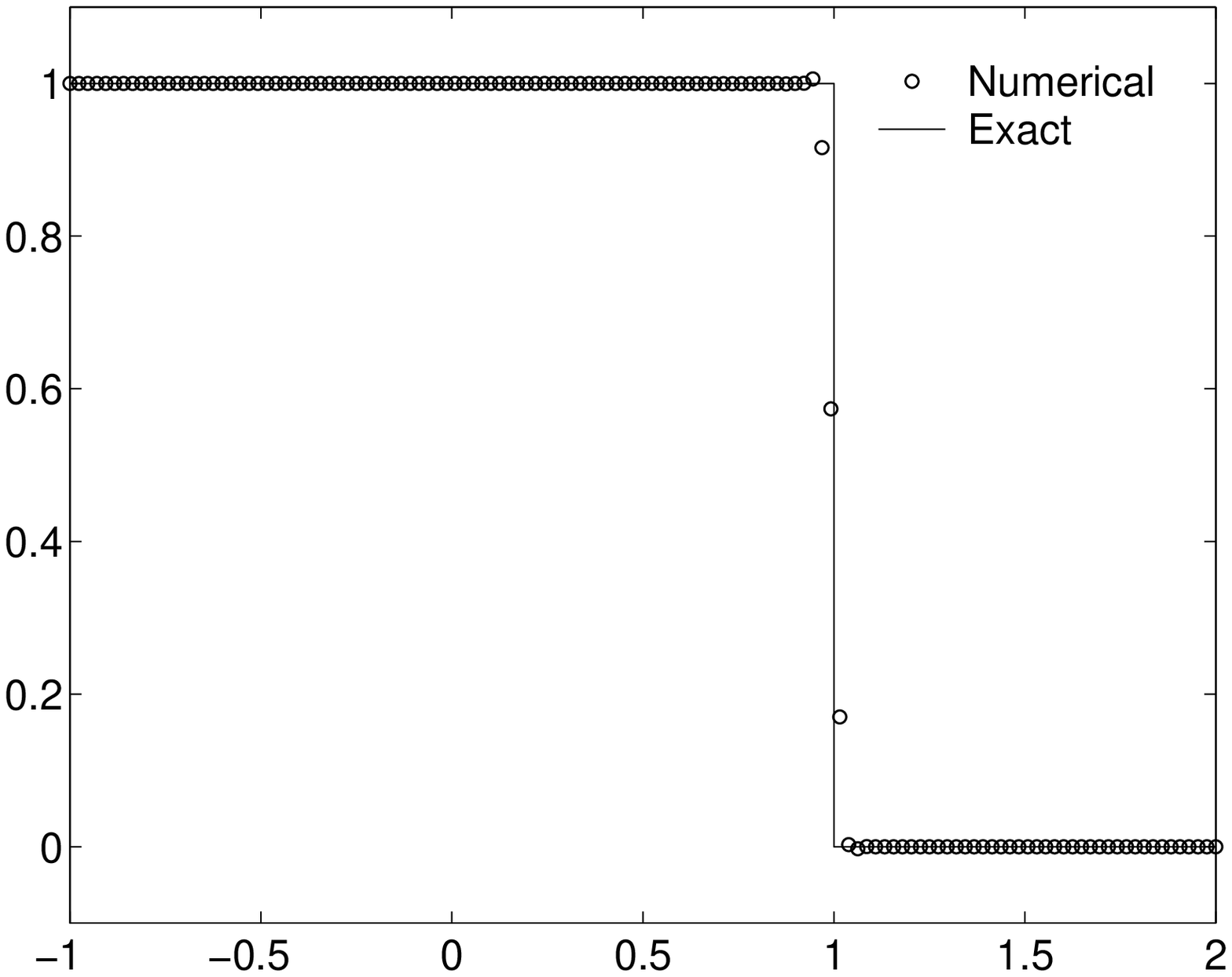}}
 \hspace{0.5cm}
 \subfigure[ ]{
     \includegraphics[width=7cm]{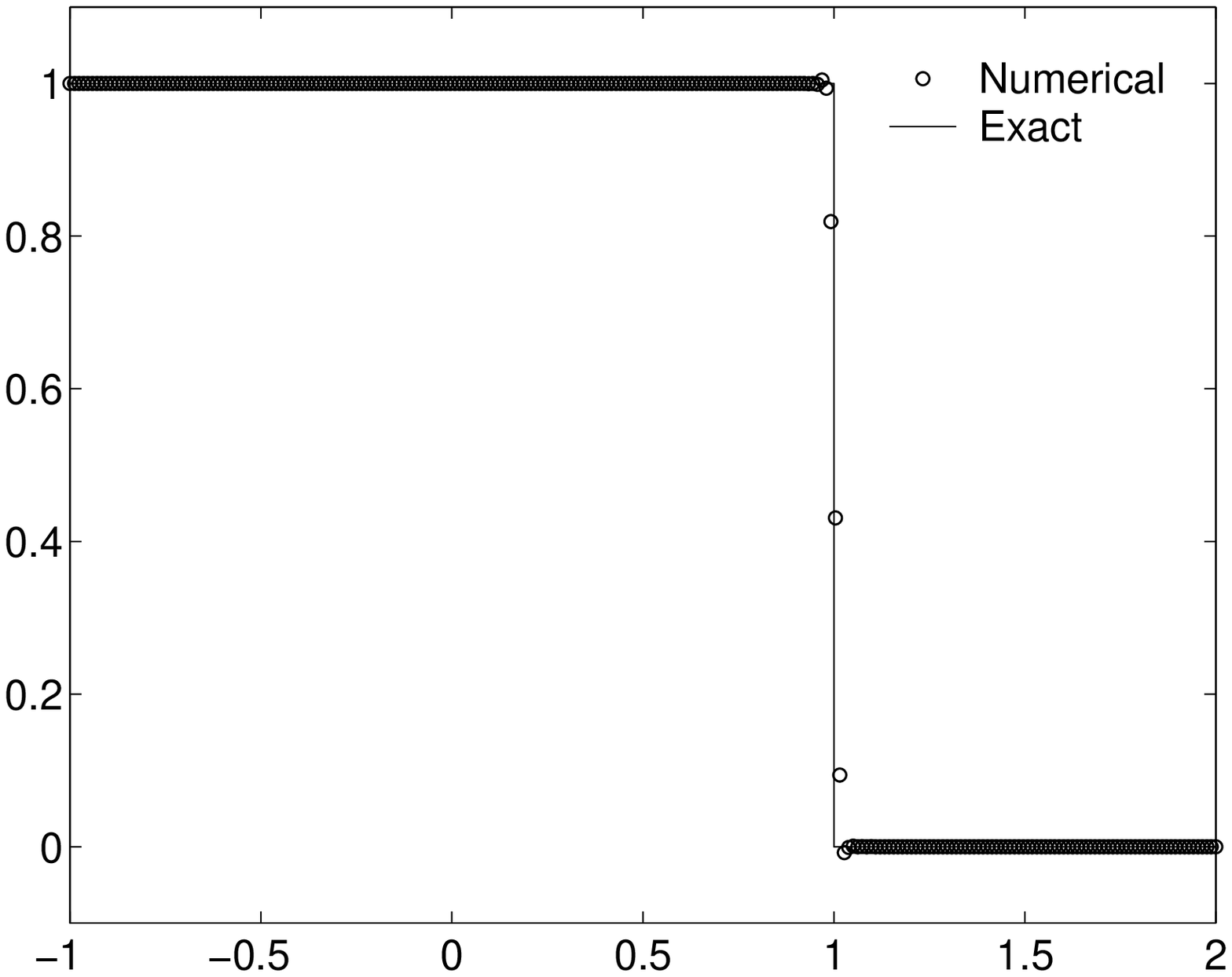}}
\caption{Inviscid Burgers' equation ($t$=2, $\Delta$ t=0.005).
        (a) 129 grid points; (b) 257 grid points.}
\label{burg1}
\end{figure}

\textit{\underline {Example 3.}} Having solved the linear equation with
different initial values, we consider the solution of the most popular
inviscid Burgers' equation, in which $f(u)=u^2/2$ and the Riemann type
initial value is given by
\begin{eqnarray}
u(x, 0)=
\l \{
   \begin{array}{ll}
   1  &  \quad  x \le 0 \\
   0  &  \quad  x>0.
   \end{array}
\r.
\end{eqnarray}
This is a standard benchmark problem in hyperbolic conservation laws
and has been considered by numerous researchers. The exact solution is
a shock wave with a constant velocity
\begin{eqnarray}
u(x, t)=
\l \{
   \begin{array}{ll}
   1  &  \quad x-St < 0 \\
   0  &  \quad x-St >0,
   \end{array}
\r.
\end{eqnarray}
where the speed of the shock front is
\begin{eqnarray}
   S=\f{1}{2}.
\end{eqnarray}
We should note that this case is a non-periodic boundary problem.
We must symmetrically double the computational domain in the $x$-direction
when we calculate the derivative of the flux. This operation is necessary
because it generates the periodic boundary condition so that the Fourier
spectral method can be used. The numerical results are plotted
in FIG. \ref{burg1} at time $t$=2.

\begin{figure}
 \centering
 \subfigure[ ]{
     \includegraphics[width=7cm]{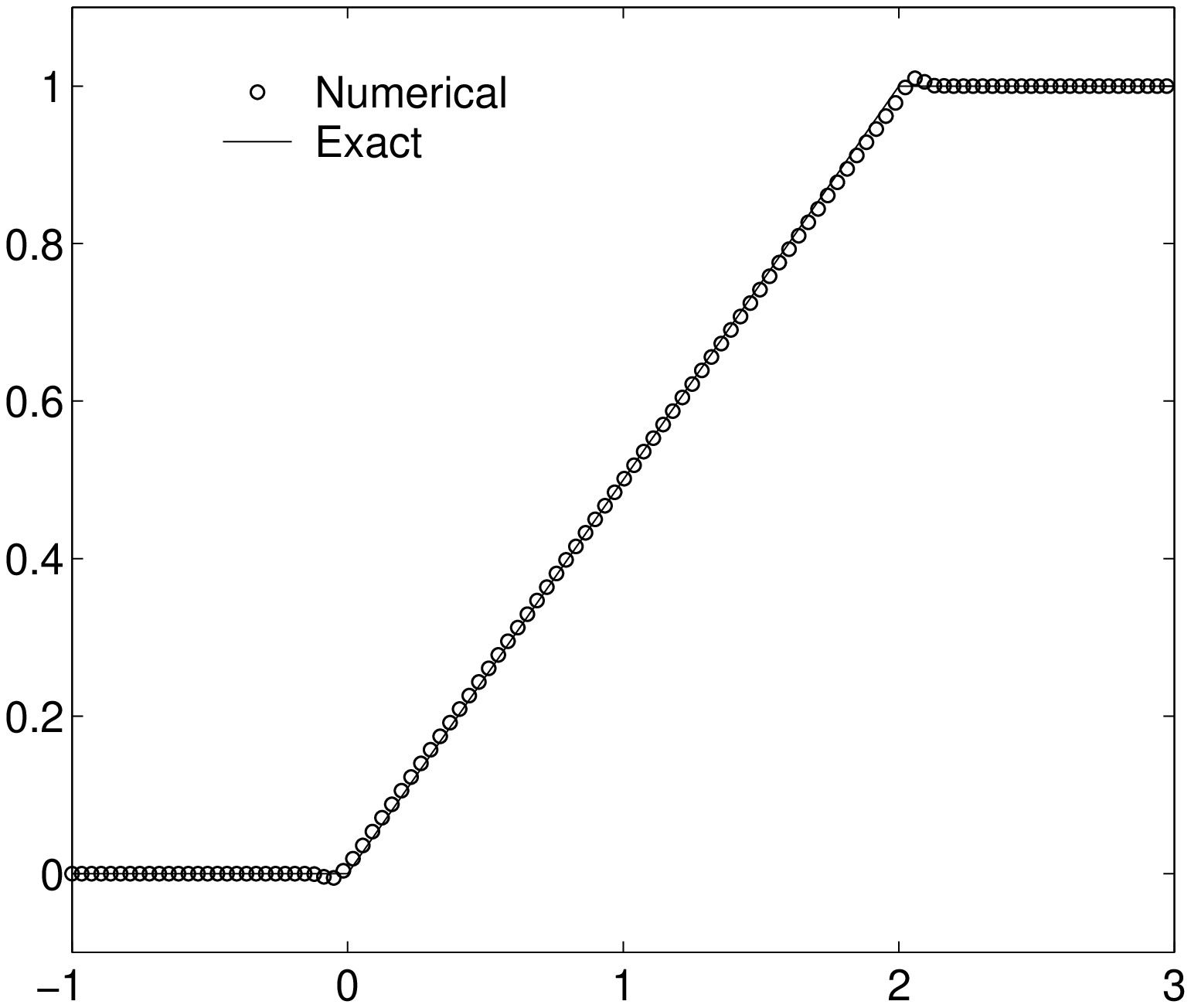}}
 \hspace{0.5cm}
 \subfigure[ ]{
     \includegraphics[width=7cm]{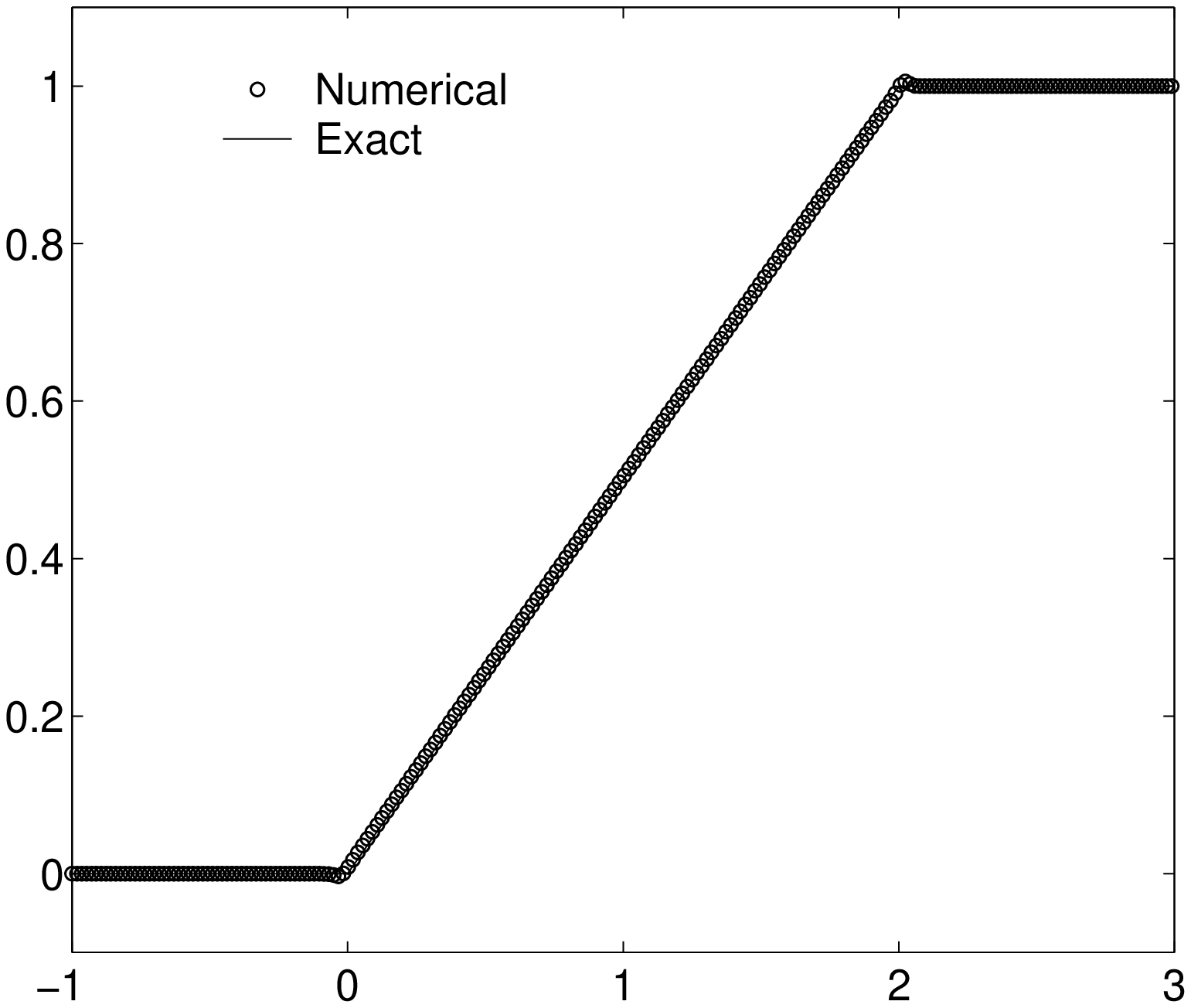}}
\caption{Inviscid Burgers' equation with rarefaction wave ($t$=2, $\Delta$t=0.005).
         (a) 129 grid points; (b) 257 grid points.}
\label{burg2}
\end{figure}

\textit{\underline {Example 4.}} Another Riemann type initial value
for the inviscid Burgers' equation with a flux of  $u^2/2$ is given by
\begin{eqnarray}
u(x,0)=
\l \{
   \begin{array}{ll}
   0  &  \quad x < 0 \\
   1  &  \quad x \ge 0.
   \end{array}
\r.
\end{eqnarray}
The exact solution of this problem is a rarefaction wave
\begin{eqnarray}
u(x, t)=
\l \{
   \begin{array}{ll}
   0         &  \quad  \f{x}{t} < 0 \\
   \f{x}{t}  &  \quad  0  <  \f{x}{t} <1 \\
   1         &  \quad  \f{x}{t} > 1.
   \end{array}
\r.
\end{eqnarray}
Numerical results are plotted in FIG. \ref{burg2} at $t$=2 with different
numbers of nodal points. Particularly, one can see that the end of the
rarefaction fan is well resolved.

\begin{figure}
 \centering
 \subfigure[ ]{
     \includegraphics[width=7cm]{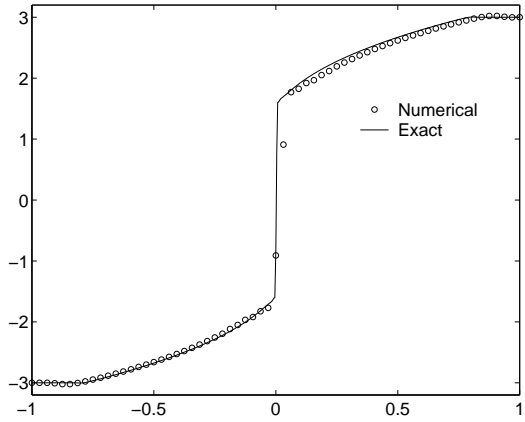}}
 \hspace{0.5cm}
 \subfigure[ ]{
     \includegraphics[width=7cm]{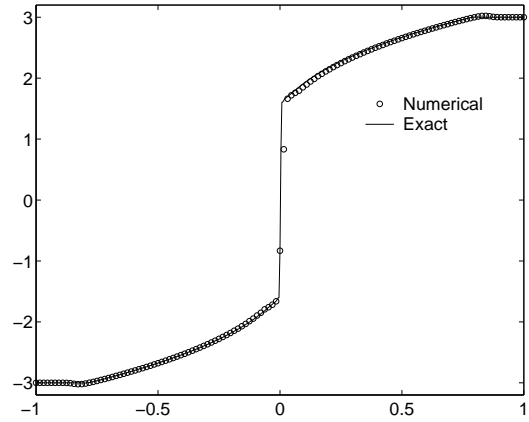}}
\caption{Inviscid Burgers' equation with non-convex flux ($t$=0.04, $\Delta$t=0.0005).
         (a) 65 grid points; (b) 129 grid points.}
\label{burg3}
\end{figure}

\textit{\underline {Example 5.}} Finally, we use non-convex flux to test the
convergence to the physically correct solution. The flux is a non-convex
function
\begin{eqnarray}
   f(u)=\f{1}{4} (u^2-1)(u^2-4)
\end{eqnarray}
with a Riemann type initial value
\begin{eqnarray}
u(x,0)=
\l \{
   \begin{array}{ll}
   -3  &  \quad x < 0 \\
    3  &  \quad x \ge 0.
   \end{array}
\r.
\end{eqnarray}
The reader is referred to \cite{harten} for more detailed information about this problem.
The results are displayed in FIG. \ref{burg3} at $t$=0.04. Again, present
spectral solver yields a satisfactory resolution.

\subsection{Euler systems in one dimension}

In this subsection, we perform numerical experiments by using the proposed
scheme for the 1D Euler equation of gas dynamics. In one dimension, the Euler
equation takes the form
\begin{eqnarray}
U_t + F(U)_x =0 \label{Euler1}
\end{eqnarray}
with
\begin{eqnarray}
U= \l ( \begin{array}{c}
         \rho   \\
         \rho u \\
              E \\
       \end{array} \r ); \quad
F(U) = \l ( \begin{array}{c}
             \rho u       \\
             \rho u^2 + p \\
             u(E+p)  \\
           \end{array} \r ),
\end{eqnarray}
where, $\rho, u, p$ and $E$ denote the density, velocity, pressure
and total energy per unit mass $E=\rho(e+(u^2)/2)$, respectively.
Here, $e$ is the specific internal energy. For an ideal gas with the constant
specific heat ratio ($\gamma =1.4$) considered here, one has
$e=p/(\gamma-1)\rho$. We consider the following two well-known
Riemann problems.

\begin{figure}
 \centering
 \subfigure[ ]{
     \includegraphics[width=7cm]{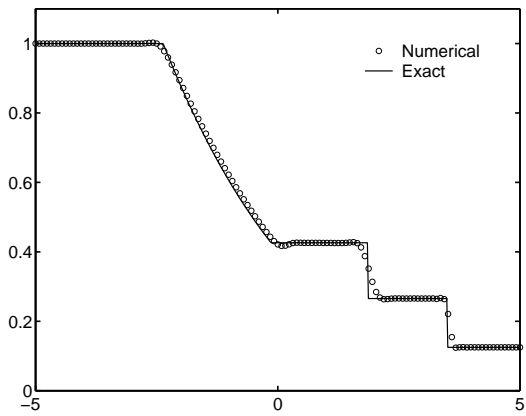}}
 \hspace{0.5cm}
 \subfigure[ ]{
     \includegraphics[width=7cm]{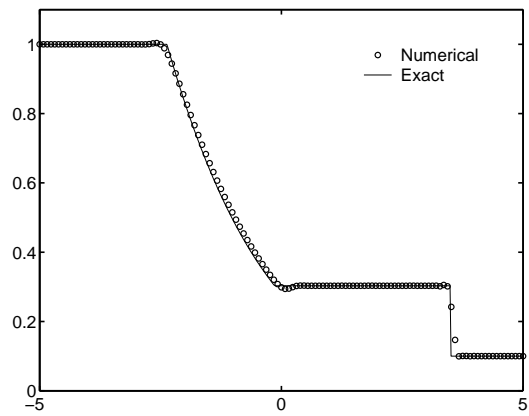}}
\caption{Sod's problem ($t$=2, $\Delta$t=0.02, 129 grid points).
         (a) Density; (b) Pressure.}
\label{sod}
\end{figure}

\begin{figure}
 \centering
 \subfigure[ ]{
     \includegraphics[width=7cm]{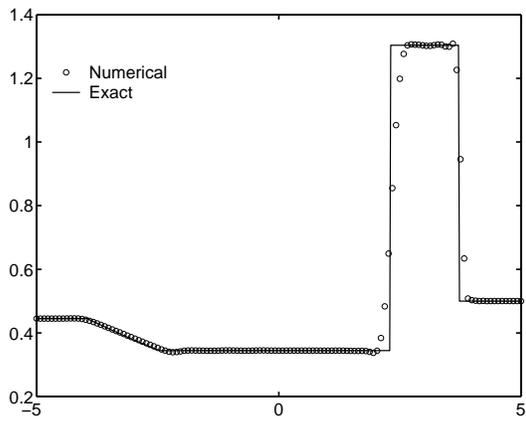}}
 \hspace{0.5cm}
 \subfigure[ ]{
     \includegraphics[width=7cm]{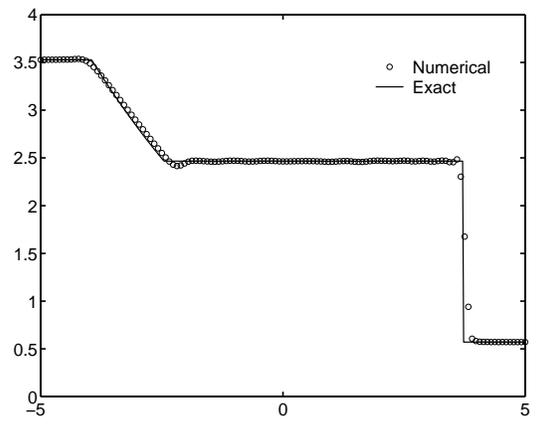}}
\caption{Lax's problem ($t$=1.5, $\Delta$t=0.02, 129 grid points).
         (a) Density; (b) Pressure.}
\label{lax}
\end{figure}

\textit{\underline {Example 6.}} We solve for the solution of two shock tube
problems with Sod's and Lax's initial conditions \cite{shu96}, i.e.,

\begin{eqnarray}
(\rho, u, p)_{t=0}= \left \{
   \begin{array}{ll}
     (1,0,1),         &   \quad  x < 0 \\
     (0.125,0,0.1),   &   \quad  \text{otherwise}
   \end{array}
   \r.
\end{eqnarray}
for the Sod problem and
\begin{eqnarray}
(\rho, u, p)_{t=0}= \left \{
   \begin{array}{lll}
    (0.445,0.698,3.528),   &   \quad x < 0 \\
    (0.5,0,0.571),         &   \quad \text{otherwise}
   \end{array}
   \right.
\end{eqnarray}
for the Lax problem. The numerical results for these two problems are shown in
FIG. \ref{sod} and FIG. \ref{lax}, respectively. We can see that the present
method gives a good resolution in both cases. Similar to all the methods depending
on the numerical dissipation, the resolution for the linear degenerate contact
discontinuity is slightly lower than that for the generic nonlinear shock wave. If
combined with Harten's artificial compression method, present spectral solver should
have better resolution for the contact discontinuity.

\textit{\underline {Example 7.}} In this example, the interaction
of an entropy wave of small amplitude with a Mach 3 right-moving
shock in a one-dimensional flow is investigated. The computational
domain is taken as $[0,9]$ and the flow field is initialized with
\begin{eqnarray} \label{inicond}
 (\rho, u, p)_{t=0}  = \left \{
         \begin{array}{lccl}
         (3.85714, & 2.629369, & 10.33333 ) &  \quad x \le 0.5 \\
         (e^{- \epsilon \mbox {sin}( \kappa x)}, & 0, & 1.0 ) &  \quad  x > 0.5 \\
        \end{array} \r.
\end{eqnarray}
where $\epsilon$ and $\kappa$ are the amplitude and the wave
number of the entropy wave before the shock. This problem is
significant due to its relevance to the interaction of
shock-turbulence. Spectral methods, by nature of their high accuracy and
negligible dispersive and dissipative errors, are ideally suited to
the numerical study of turbulence and have already been widely applied.
Here we would like to argue the feasibility of
the present spectral method to the shock-turbulence interaction. In our test, we vary
the wavenumber of the pre-shock wave while keeping its amplitude unchanged
($\epsilon=0.01$). As the wavenumber increases, the problem becomes more and more
challenging because the amplified high-frequency entropy waves after the shock are
always mixed with spurious oscillations. As the wavenumber increases, it
is very difficult to distinguish the high-frequency waves from the spurious
oscillations. A low-order scheme may dramatically damp the transmitted
high frequency waves. Even some popular high-order schemes also
encounter the difficulty in preserving the amplitudes of the
entropy waves due to their excessive dissipation with a given mesh size.
Therefore, a successful shock-capturing method should be able to
eliminate Gibbs' oscillation while capturing the shock and preserving
the high-frequency entropy waves.

In our test, we let the shock move from $x=0.5$ to
$x=8.5$. For the purpose of comparison with the previous results,
we only give the results at the interval $[4.0, 9.0]$. Also, in order
to discharge the transient waves due to the non-numerical
initial shock profile, we plot the length of the amplified
entropy waves in the same manner as that in Ref. \cite{shu96}. Furthermore,
we must point out two nontrivial remarks for the following cases
with different $\kappa$. One is that we have replotted the final results
in a denser grid, otherwise the bounded strip will not be fully spanned
due to the fact that no enough grid values locate near all the
peaks /valleys of the wave. Another is that we have used a local filter
in order to eliminate the oscillations near the shock when we
present final results after the computation.

\begin{figure}
 \centering
 \subfigure[$\kappa=13$, $N=513$]{
     \includegraphics[width=12cm]{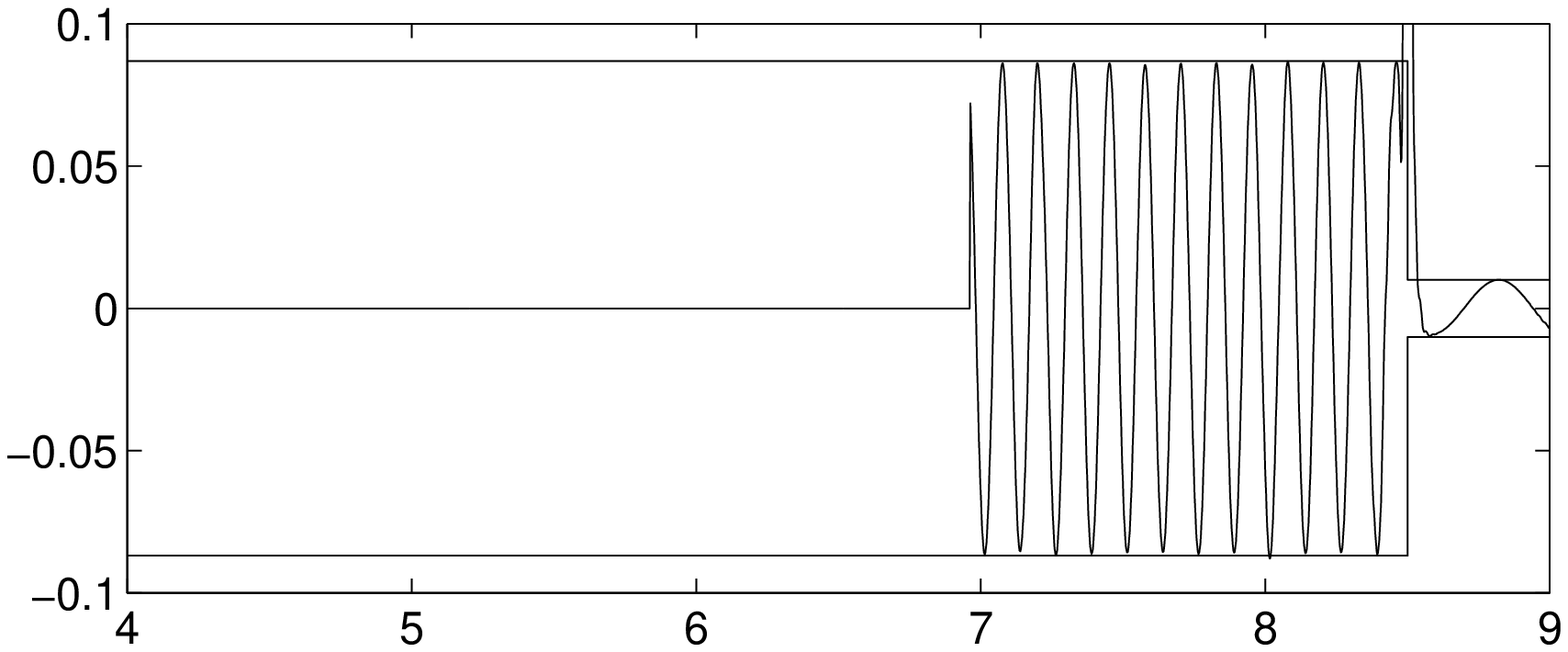}} \\

 \subfigure[$\kappa=26$, $N=1025$]{
     \includegraphics[width=12cm]{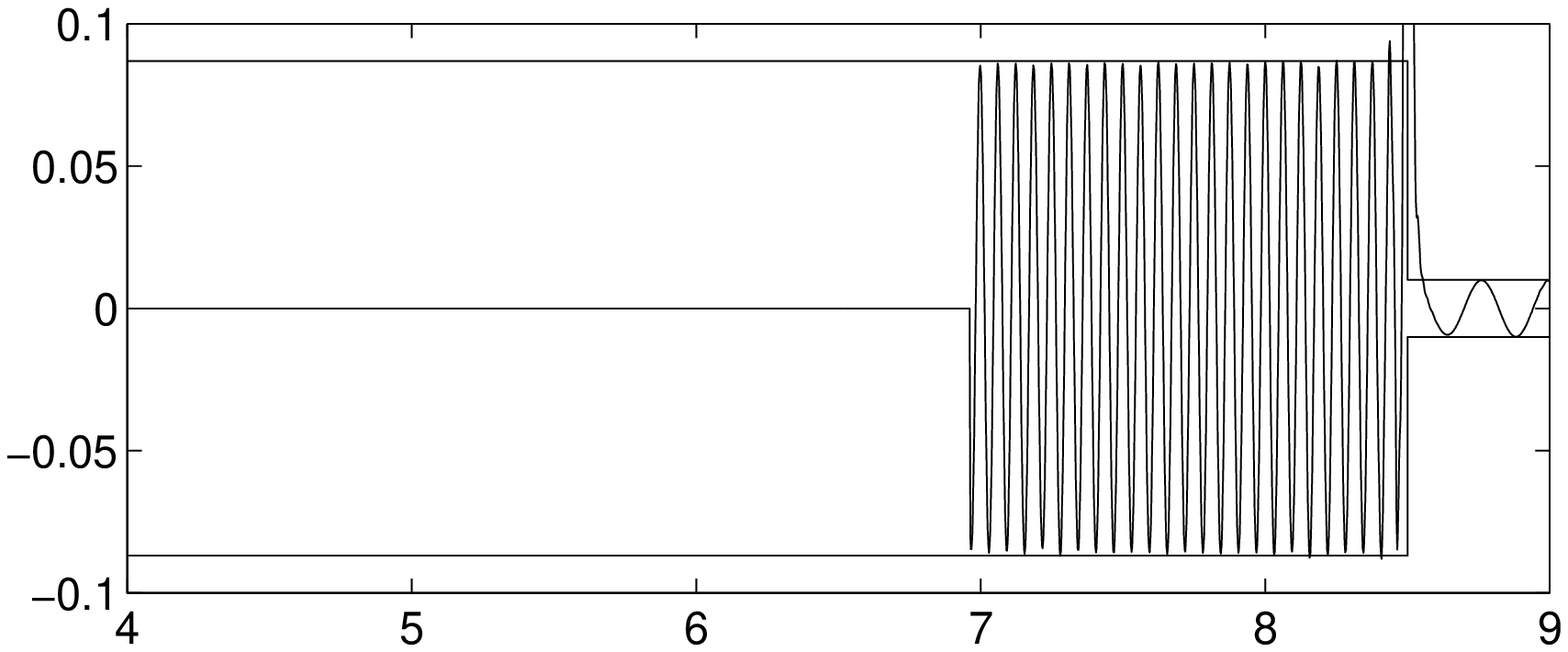}} \\

 \subfigure[$\kappa=39$, $N=2049$]{
     \includegraphics[width=12cm]{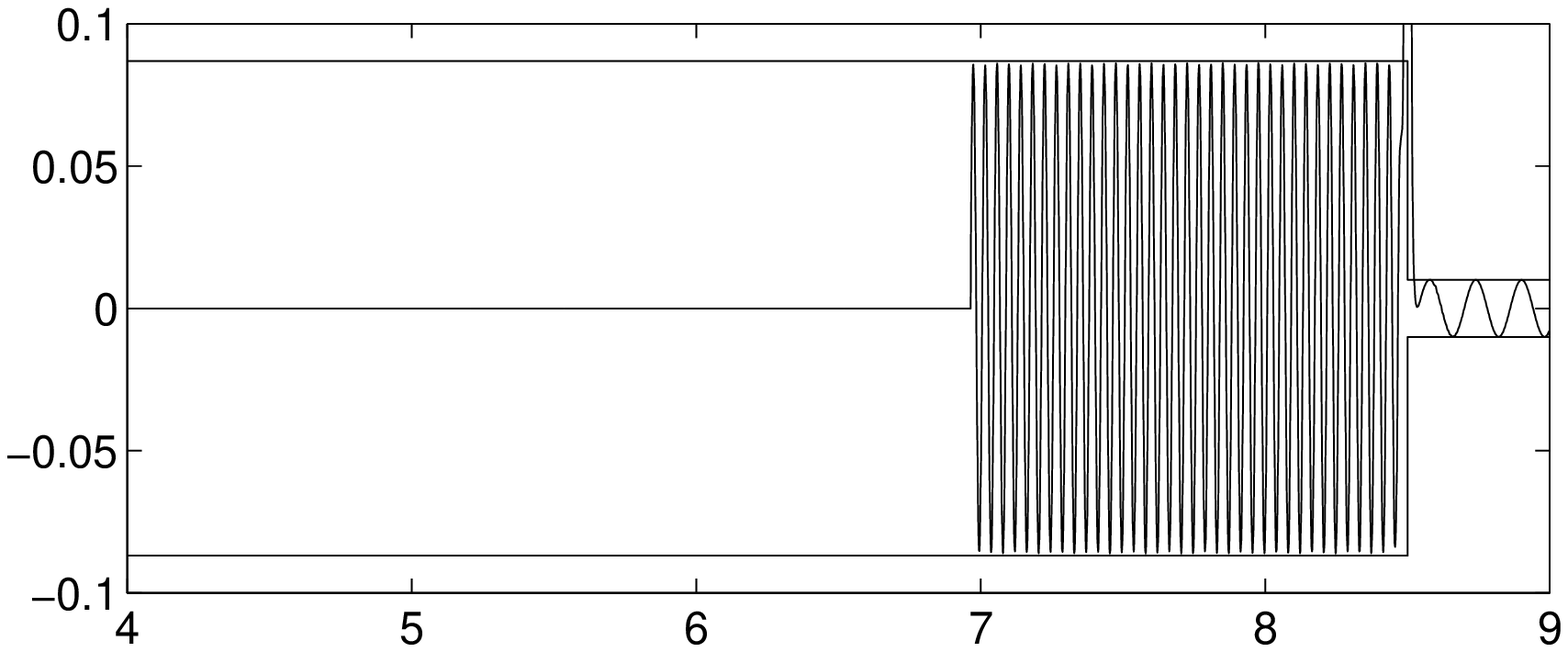}} \\

 \subfigure[$\kappa=52$, $N=2049$]{
     \includegraphics[width=12cm]{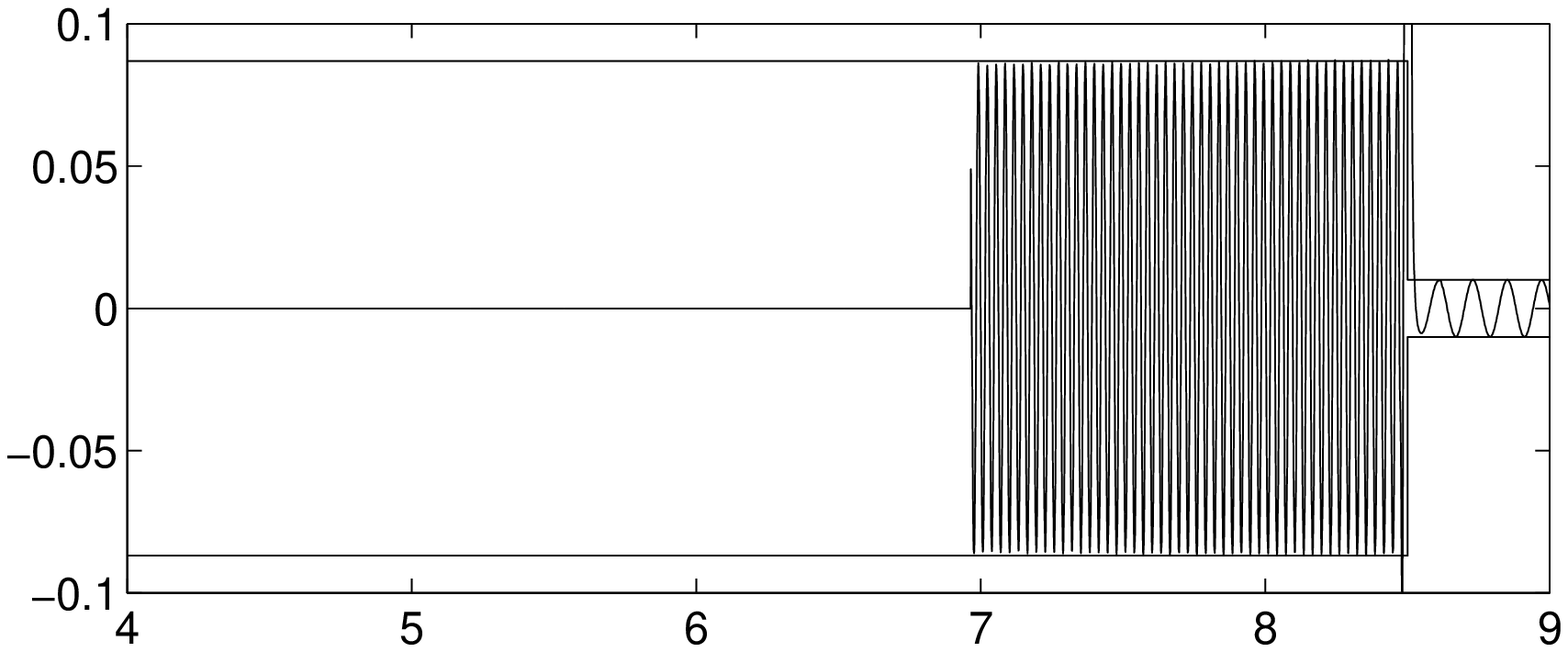}} \\
\caption{1D shock entropy wave interaction. Amplitude of entropy waves.}
\label{shockentropy}
\end{figure}

First, we consider $\kappa=13$. In this case, the computational
domain is deployed with 513 grids points. Such a mesh is suitable
for the Fourier spectral method.
The numerical results of the amplitude of the entropy waves
are shown in FIG. \ref{shockentropy}(a). It can
be seen that the shock and the generated entropy waves are ideally
captured. The entropy waves fully span the strip bounded by
two solid lines at $\pm$ 0.08690716, which show the amplitude of
the amplified entropy waves predicted by the linear analysis.
Since post-shock waves are monochromic, a simple calculation
indicates the point per wavelength (PPW) value of these waves
is 7.125. Such value is larger than the spectral resolution (PPW=2.0)
but it is much smaller than that obtained by most shock-capturing
schemes. The transmitted waves suffer the similar smearing as
the contact waves. The shock wave, on the contrary, is self-repairable due
to its compression nature. This observation verified that the numerical
dissipation needed for shock-capturing will degrade the resolution
of the spectral method. It also verified that this degradation can
be isolated and alleviated by the DSC-RSK filter.

We next double the wave number $\kappa$ up to 26. A mesh of 513
grid points is not enough for simultaneous shock capturing and
high-frequency wave resolving. Therefore, we increase the mesh size to
$N=1025$ for this computation. The results for this case
are presented in FIG. \ref{shockentropy}(b). It is seen that
no obvious deterioration happens to
the results and the resolution in the post-shock waves is still
satisfactory. The fullness of the wave is comparable with that
of $\kappa=13$. But the first two post-shock waves have been
slightly polluted.

When the wave number $\kappa$ is increased to 39 and 52, we
use 2049 grid points to resolve such high-frequency waves.
FIG. \ref{shockentropy}(c) and FIG. \ref{shockentropy}(d)
show the results of $\kappa=39$ and $\kappa=52$, respectively.
The same PPW value $7.125$ is maintained in this high-frequency
case. The good resolution of the present method in simulating
this 1D shock entropy interaction indicates that our scheme is
capable of and efficient in distinguishing the high-frequency entropy
waves from  spurious oscillations.

\begin{figure}
 \centering
 \subfigure[ ]{
     \includegraphics[width=7cm]{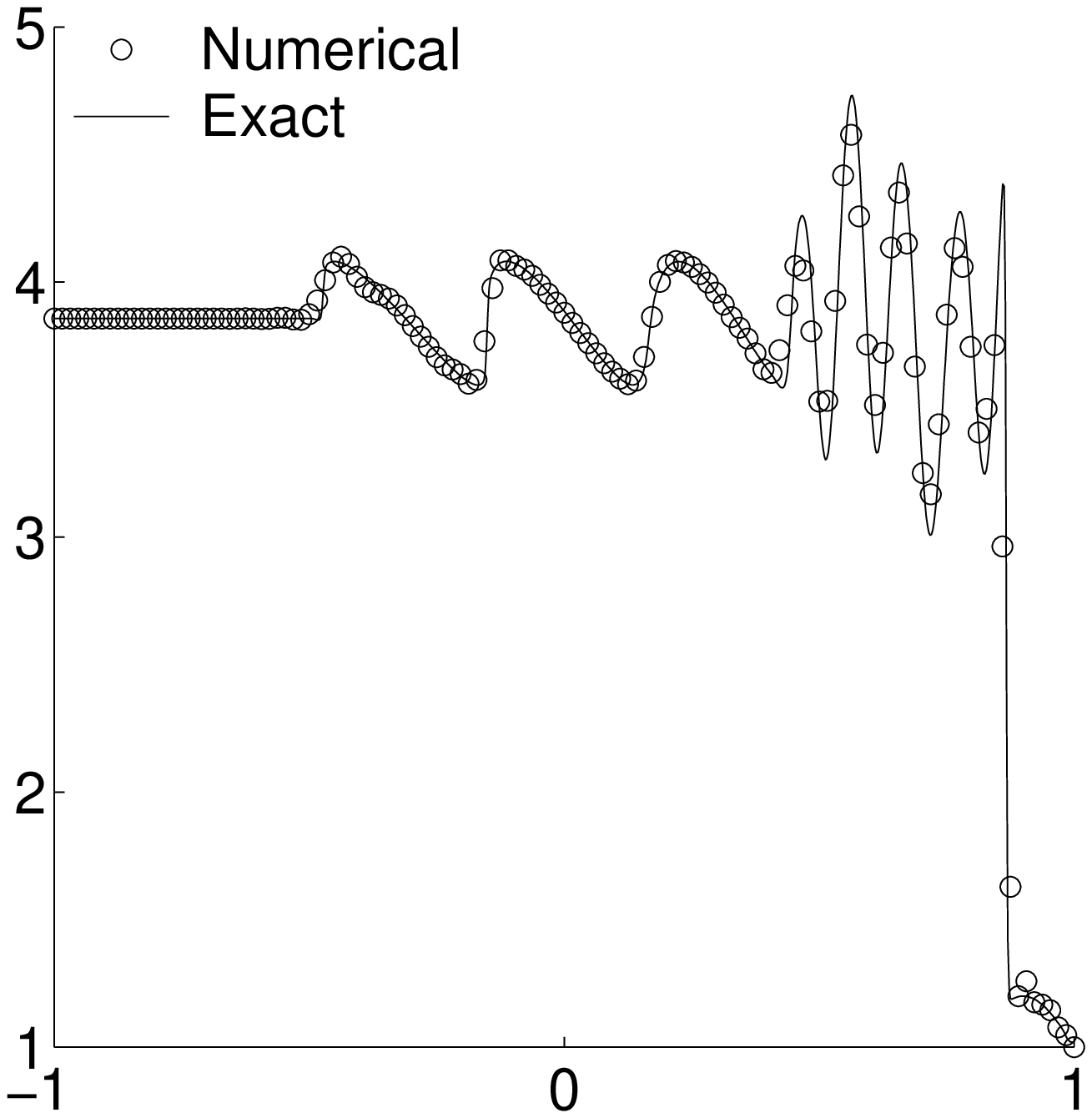}}
 \hspace{0.1cm}
 \subfigure[ ]{
     \includegraphics[width=7cm]{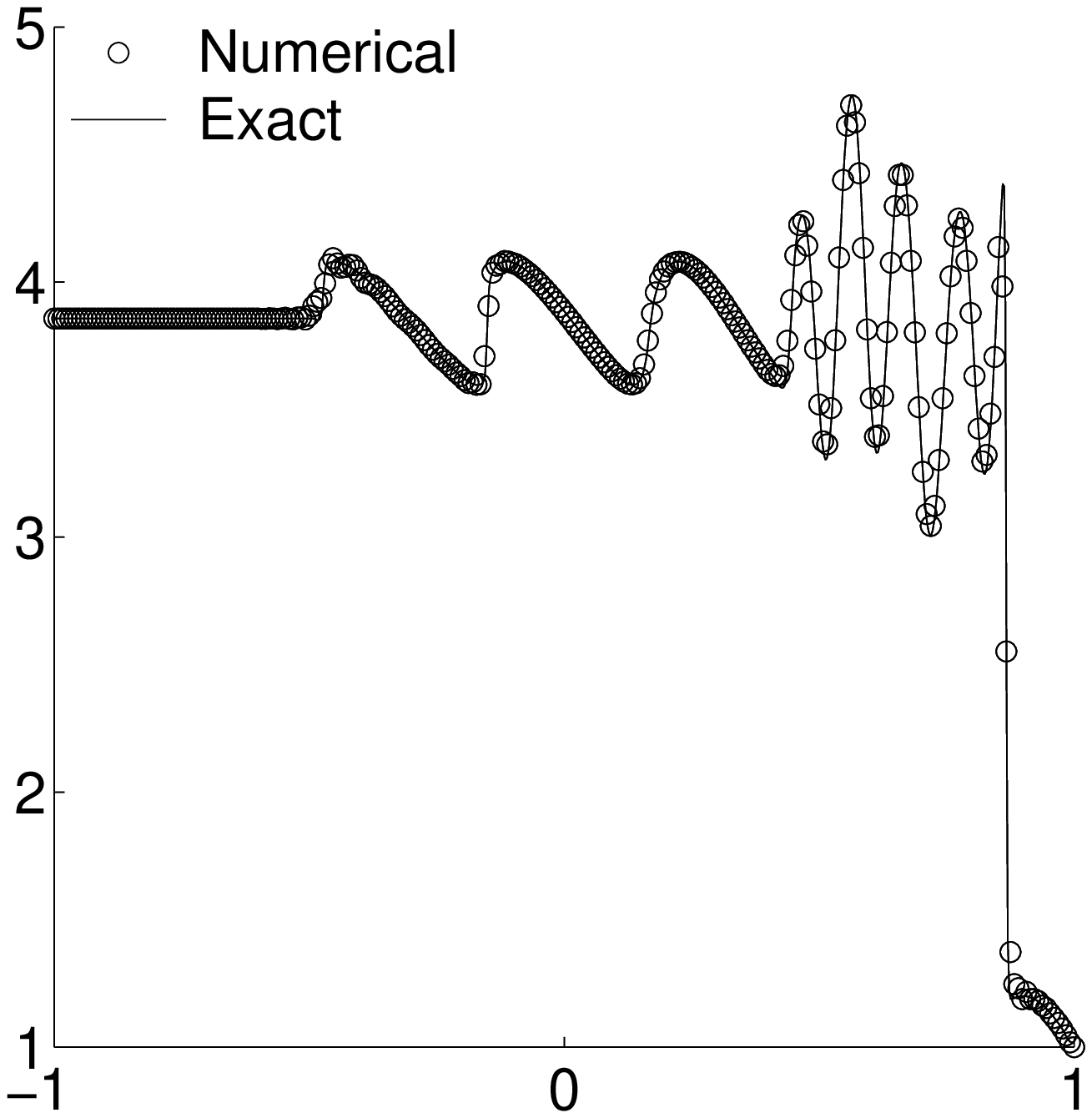}}
\caption{Shu-Osher problem ($t$=0.47).
         (a) 128 cells; (b) 256 cells.}
\label{shu}
\end{figure}

\textit{\underline {Example 8.}} Results are now shown for the problem
of Shu and Osher which also describes the interaction of an entropy
sine wave with a Mach 3 right-moving shock. The computational domain
is taken as $[-1,1]$ and the flow field is initialized with
\begin{eqnarray} \label{inicond1}
(\rho, u, p)_{t=0}  =
\left \{
     \begin{array}{lcrr}
         (3.85714, & 2.629369, & 10.33333 ) &  \quad x \le -0.8\\
         (1.0+\epsilon \mbox {sin}( \kappa \pi x), & 0, & 1.0 ) &  \quad  x > -0.8,
     \end{array}
 \r.
\end{eqnarray}
where $\epsilon$=0.2 and $\kappa$=5. FIG. \ref{shu}(a) and
FIG. \ref{shu}(b) show the solution of the density with 128 and 256
cells, respectively. The `Exact' solution here is the solution
computed by the third-order ENO with 1200 cells. As evident in
the solution, the complicated flow field behind the shock is well
resolved and the shock remains sharp. However, small oscillations
near the shock still exist due to the fact that we do not
attempt to locate this shock intentionally. If we predict the
position of the shock and apply a local filter near the shock at
the final time step, it is believed that more satisfactory results
can be obtained.

\subsection{Euler systems in two dimensions}

In two dimensions, the Euler equation for gas dynamics in
a vector notation takes the conservation form
\begin{eqnarray}
U_t + F(U)_x + G(U)_y=0 \label{Euler}
\end{eqnarray}
with
\begin{eqnarray}
U= \l ( \begin{array}{c}
         \rho   \\
         \rho u \\
         \rho v \\
              E \\
       \end{array} \r ); \quad
F(U) = \l ( \begin{array}{c}
             \rho u       \\
             \rho u^2 + p \\
             \rho u v     \\
             u(E+p)  \\
           \end{array} \r ); \quad
G(U) = \l ( \begin{array}{c}
             \rho v       \\
             \rho u v     \\
             \rho v^2 + p \\
             v (E+p)  \\
           \end{array} \r ),
\end{eqnarray}
\begin{eqnarray}
p=(\gamma-1)(E-\f{1}{2} \rho(u^2+v^2)),
\end{eqnarray}
where $(u, v)$ is the 2D fluid velocity and $\gamma=1.4$ is used in
all computations.

\textit{\underline {Example 9.}} As the first example, we use our
spectral scheme to study the interaction between a normal shock and a
weak entropy wave which makes an angle $\theta \in (0, \pi/2)$
against the $x$-axis. If $\theta=0$, we have essentially
the 1D problem (see example 7). The initial conditions are
defined as follows: the right state of the shock
is given as $(\rho_r, u_r, v_r, p_r)$=(1, 0, 0, 1) with a given shock Mach
number $M_s=3$. We add a small entropy wave to the flow on the right
of the shock which is equivalent to changing only the density of the
flow on the right of the shock. In our problem we change the density
$\rho_r$ in the right state of the shock by multiplying $\rho_\epsilon$
\begin{eqnarray}
   \begin{array}{ll}
       \rho=           &  \rho_r\rho_\epsilon\\
       \rho_\epsilon=  &  e^{-(\epsilon/p_r) \sin
                         (\kappa(x \cos\theta+y \sin \theta))}
   \end{array}
\end{eqnarray}
with $\epsilon$ and $\kappa$ being the amplitude and wavenumber, respectively.
In order to carry out a long time integration and to enforce
periodic boundary condition in $y$-direction, the computational
domain is taken to be  $[0,9] \times [0,\f{2\pi}{\kappa \sin \theta}]$.
We initially position the normal shock at $x=0.5$ and allow
it to move up to $x=8.5$. In our simulation, we adopt
$\epsilon=0.1$, $\kappa=15$, and $\theta=30^\circ$.
The performance is measured by drawing the maximum amplitudes of the
amplified entropy waves in the $y$-direction for all fixed grid values
$x \in [7.4, 8.4]$, and comparing them with the amplitude predicted by
the linear analysis, which is 0.08744786.

\begin{figure}
\centering
    \includegraphics[width=10cm]{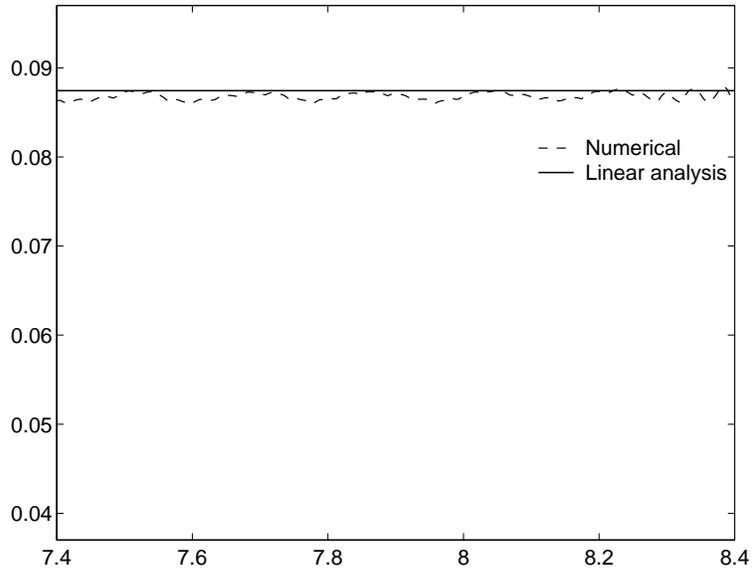}
\caption{2D shock entropy wave interaction ($N=513$). Amplitude of entropy waves.}
\label{enshock}
\end{figure}

We use 513 grid points in the $x$-direction and 32 grid points in
the $y$-direction in our test. In order to avoid peak/valley losses
of the entropy waves due to insufficient grid points near these positions, we
interpolate our final results to a denser grid when plotting them. FIG. \ref{enshock}
displays the amplitudes of the amplified entropy waves.
Obviously, the compressed entropy waves are well resolved
in the post-shock regime even though small oscillations
remain near the location of the shock. It is believed that such a
trivial flaw is acceptable.

This example has further demonstrated the capability of our spectral
scheme for capturing small scale waves in the presence of a shock.

\textit{\underline {Example 10.}} This model problem describes the
interaction between a stationary shock and a vortex. It has many
potential applications and has been treated by many researchers using
various techniques. Early numerical studies in this area using
spectral methods focused primarily on shock-fitting techniques.
However, for the cases in which a vortex causes a shock to
deform at the point of bifurcation, it is generally difficult for
a shock-fitting method to apply. This problem was successfully simulated
in Refs. \cite{cai93,wai94} by using the
Chebyshev spectral method. In this paper, we also adopt this model to
access the capability of the proposed spectral shock-capturing method.
This problem is defined as follows. The initialization is imposed on
the domain $[0, 2] \times [0, 1]$ with a stationary normal shock
at $x=0.5$ and a Mach 1.1 flow at the inlet. The right state of
the shock is given as $(\rho_r, u_r, v_r, p_r)=(1, 1.1\sqrt{\gamma}, 0, 1)$.
A vortex is generated by a perturbation to the velocity $(u, v)$,
temperature $T$ and entropy $S$ which is centered at $(x_c, y_c)=(0.25, 0.5)$
and in the forms
\begin{eqnarray}
u' & = &  \epsilon \tau e^{\alpha(1-\tau^2)} \sin \theta \\
v' & = & -\epsilon \tau e^{\alpha(1-\tau^2)} \cos \theta \\
T' & = & -\f{(\gamma-1) \epsilon^2 e^{2 \alpha(1-\tau^2)}}{4\alpha \gamma} \\
S' & = &  0
\end{eqnarray}
where, $\tau=\f{r}{r_c}$, $r=\sqrt{(x-x_c)^2+(y-y_c)^2}$
and $\theta=
\mbox {tan}^{-1}((y-y_c)/(x-x_c))$.
Here $\epsilon$ denotes the strength of the vortex,
$\alpha$ is the decay rate of the vortex. Note that we can obtain
perturbations in $\rho$ and $p$ according to the relations of
$T=p/\rho$ and $S={\rm ln}\f{p}{\rho^{\gamma}}$. In our numerical experiment,
we choose $\epsilon=0.3$, $r_c=0.05$ and $\alpha=0.204$.

\begin{figure}
    \centering
    \subfigure[$t$=0.05, 30 contours]{
        \includegraphics[width=5cm]{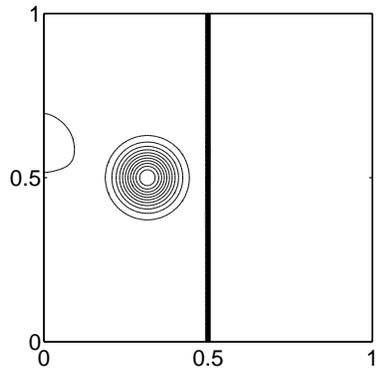}}
    \hspace{1cm}
    \subfigure[$t$=0.20, 30 contours]{
        \includegraphics[width=5cm]{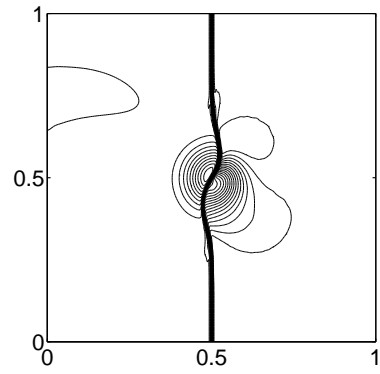}}

    \subfigure[$t$=0.35, 30 contours]{
        \includegraphics[width=5cm]{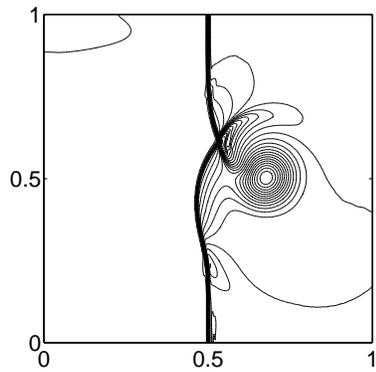}}
    \hspace{1cm}
    \subfigure[$t$=0.60, 90 contours from 1.01 to 1.37]{
        \includegraphics[width=5cm]{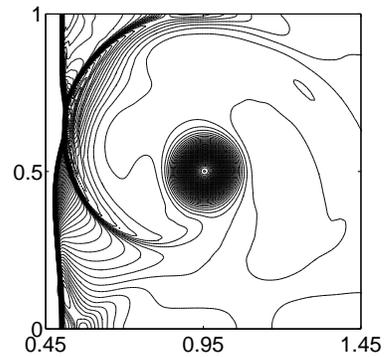}}

    \subfigure[$t$=0.80, 30 contours from 1.01 to 1.29]{
       \includegraphics[width=11cm]{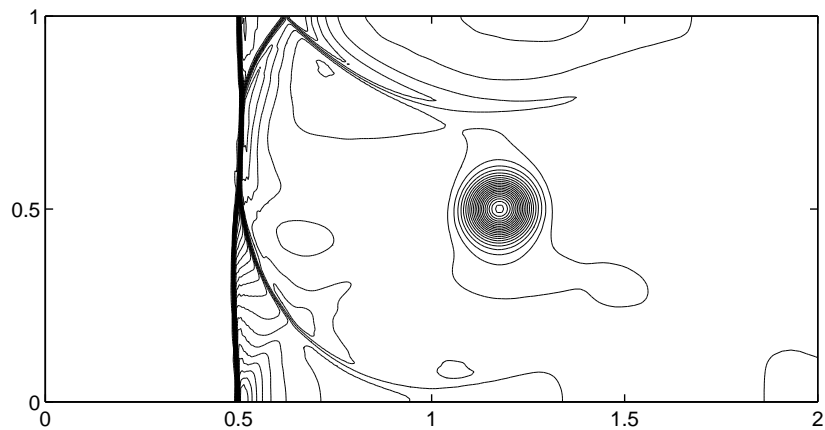}}

    \caption{2D shock vortex interaction (Mach=1.1, $\epsilon=0.3$, CFL=0.5).
             Pressure profiles. }
    \label{vorshock}
\end{figure}

A reflective boundary is imposed on both the upper and lower
boundaries. We solve the Euler equations on a 257 $\times$ 129
Cartesian grid, which is uniform in $y$-direction and refined
in $x$-direction around the shock using the Roberts
transformation \cite{roberts}. Once again, we emphasize that
the DSC-RSK lowpass filter is used in the physical domain
directly for this problem due to the presence of the reflective boundaries.
FIGs. \ref{vorshock}(a-e) show a time evolution of shock-vortex
interaction in terms of the pressure contours. We can see that even
at $t$=0.8 the solution (FIG. \ref{vorshock}(e)) is almost free of
oscillations and the fine scale features are captured very well.
Particularly, the shock bifurcation reflected on the top boundary
is shown clearly.

\begin{figure}
\centering
    \includegraphics[width=10cm]{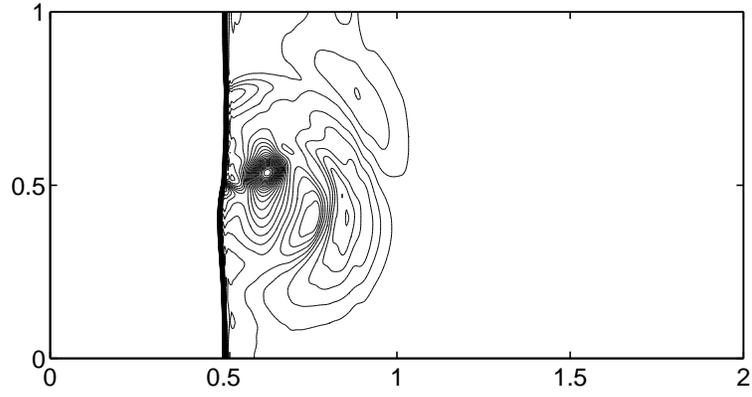}
\caption{2D shock vortex interaction ($t$=0.2, Mach=3.0, $\epsilon=0.6$, CFL=0.5).
         Pressure profiles.}
\label{vorshock3}
\end{figure}

In order to explore the potential of our algorithm further, we increase
the flow Mach number up to 3.0 and the strength of the vortex
$\epsilon$ up to 0.6. The results at $t$=0.2 are shown in
FIG. \ref{vorshock3}. Deformation and bifurcation of
the shock are observed clearly. These results indicate that the
applicability of our spectral shock-capturing scheme is not
limited to low Mach number and low vortex strength. Flows of high Mach
number and high vortex strength can also be successfully simulated.

\textit{\underline {Example 11.}} In this case, we investigate the
advection of an isentropic vortex in a free stream. As the exact
solution of this problem is available, it is a good test for
examining the accuracy and stability of shock capturing schemes.

Consider a mean flow of $(\rho_{\infty},u_{\infty},v_{\infty},
p_{\infty},T_{\infty})=(1,1,1,1,1)$ with a periodic boundary
condition in both directions. At $t_0$, the flow is
perturbed by an isentropic vortex $(u' ,v' ,T')$
centered at $(x_0,y_0)$, having the forms of
\begin{eqnarray}
u' & = & -\f{\lambda}{2 \pi}(y-y_0) e^{\eta(1-s^2)}, \\
v' & = & ~~\f{\lambda}{2 \pi}(x-x_0) e^{\eta(1-s^2)}, \\
T' & = & -\f{(\gamma-1) \lambda^2}{16 \eta \gamma \pi^2} e^{2 \eta (1-s^2)}.
\end{eqnarray}
Here, $s=\sqrt{(x-x_0)^2+(y-y_0)^2}$ is the distance to the
vortex center; $\lambda$ is the strength of the vortex and
$\eta$ is a parameter determining the gradient of the solution.
In our test, we choose them as $\lambda=5$ and $\eta=1$ unless we
further specify. Note that for an isentropic flow,
relations $p=\rho^{\gamma}$ and $T=p/\rho$ are valid.
Therefore, the perturbed $\rho$ is required to be
\begin{eqnarray}
\rho=(T_{\infty} + T')^{1/(\gamma-1)} =
\l [ 1-\f{(\gamma-1) \lambda^2}{16 \eta \gamma \pi^2} e^{2 \eta (1-s^2)}
\r ] ^{1/(\gamma-1)}.
\label{density}
\end{eqnarray}

The computational domain is taken to be $[0, 10]\times [0, 10]$ with the
center of the vortex being initially located at $ (x_0, y_0)=(5, 5)$, the
geometrical center of the computational domain.

Since there is no presence of shock in this problem, spectral method
on its own can already provide excellent results if the time integration
period is small enough. As such, the lowpass filter is not needed
in an initial time period. However, as time progresses,
errors would accumulate rapidly and the computation could become
unstable if the lowpass filter is not activated to efficiently control the
dramatically nonlinear growth of the errors. In our numerical
experiment which is not shown here, the computation blows up at around
$t$=13 if no filter is used.

Two experiments are conducted in this study. One is to examine the accuracy
of our spectral code and  the other is to investigate the stability of our
scheme for a long time integration. For the first experiment, we compute the
solution of $\rho$ at the time $t$=2 using three sets of meshes
($N=32, 64, 128$).  To make the spatial error dominant, we optimize
the CFL number to be 0.01.

Two error measures, $L_1$ and $L_2$, have been used. To be
consistent with the literature \cite{eric01}, two errors used in
this problem are defined as
\begin{eqnarray}
L_1 & = & {1\over (N+1)^2}
  \sum_{i=0}^{N} \sum_{j=0}^{N} | f_{i, j}- \bar{f}_{i, j} |  \\
L_2 & = & {1\over (N+1)} \sqrt{\sum_{i=0}^{N} \sum_{j=0}^{N}
               | f_{i,j}- \bar{f}_{i,j} |^2},
\end{eqnarray}
where $f$ is the numerical result and $\bar{f}$ the exact
solution. Note that they are not the standard definitions.
The errors for the density with respect to the exact
solution are listed in Table I in which we also list the results
of some other schemes for a comparison. As the table highlights,
the accuracy of the Fourier spectral method is comparative to
that of the DSC method \cite{yangsy}, while both Fourier spectral and DSC
methods are much more accurate than the fourth-order accurate, conservative
central scheme (C4).

Our next numerical experiment concerns the performance of our
spectral code for a long time integration. Here, we compute the
solution of $\rho$ at the time $t$=100, 200, 400, 600, 800 and 1000,
with the grid $N=64$. Note that we choose the gradient of
the solution $\eta=0.5$ in this case. We list the $L_1$ and $L_2$
errors for the density in Table II. It can be seen that our spectral
code is highly accurate even after such a long time evolution.
The results indicate that our spectral scheme together with
the DSC lowpass filter is reliable for the numerical approximation
of conservation laws.

\textit{\underline {Example 12.}} 
Finally, we consider the problem
of a supersonic flow past a cylinder \cite{shu96}. This example is used to
mainly examine the ability of our DSC filter spectral scheme
in handling a non-rectangular domain.

The physical domain is on the $x-y$ plane, and the computational
domain is chosen to be [0, 1] $\times$ [0, 1] on the $\xi-\eta$ plane.
The mapping between the computational domain and the  physical domain is
\begin{eqnarray}
   \begin{array}{lr}
     x=  &  -(R_x-(R_x-1)\xi)\mbox {cos}(\theta(2\eta-1)) \\
     y=  &   (R_y-(R_y-1)\xi)\mbox {sin}(\theta(2\eta-1))
    \end{array}
\label{eq:thansform}
\end{eqnarray}
where we take $R_x=3$, $R_y=6$ and $\theta=\frac{5\pi}{12}$.

\begin{figure}
 \centering
 \subfigure[ ]{
     \includegraphics[width=4cm]{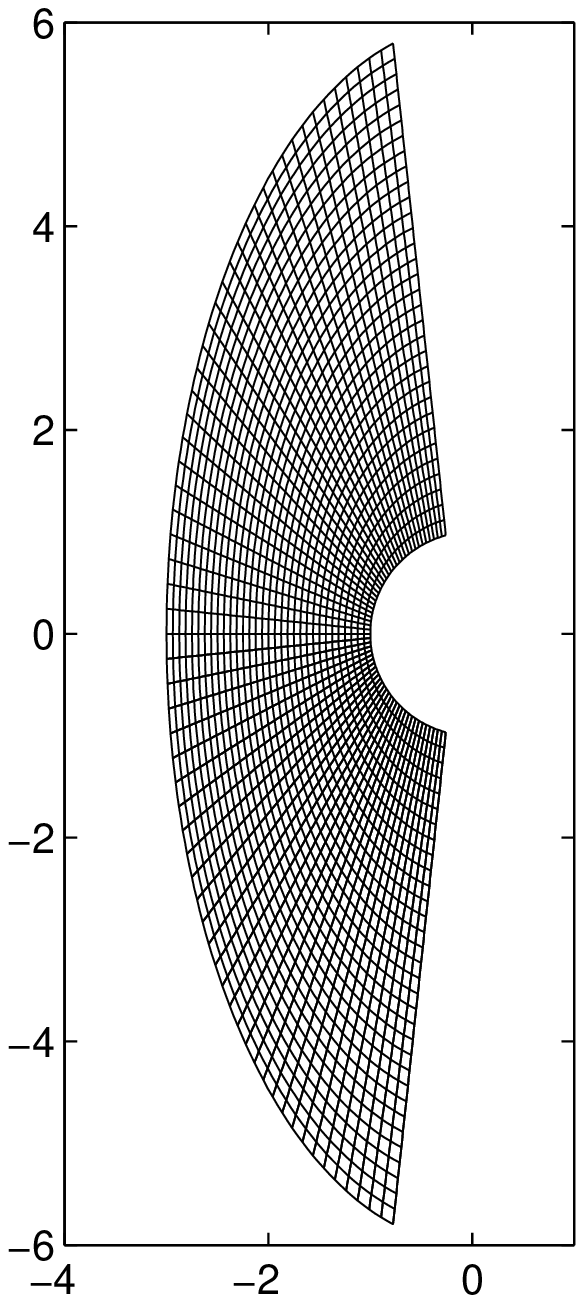}}
 \hspace{0.4cm}
 \subfigure[ ]{
     \includegraphics[width=4cm]{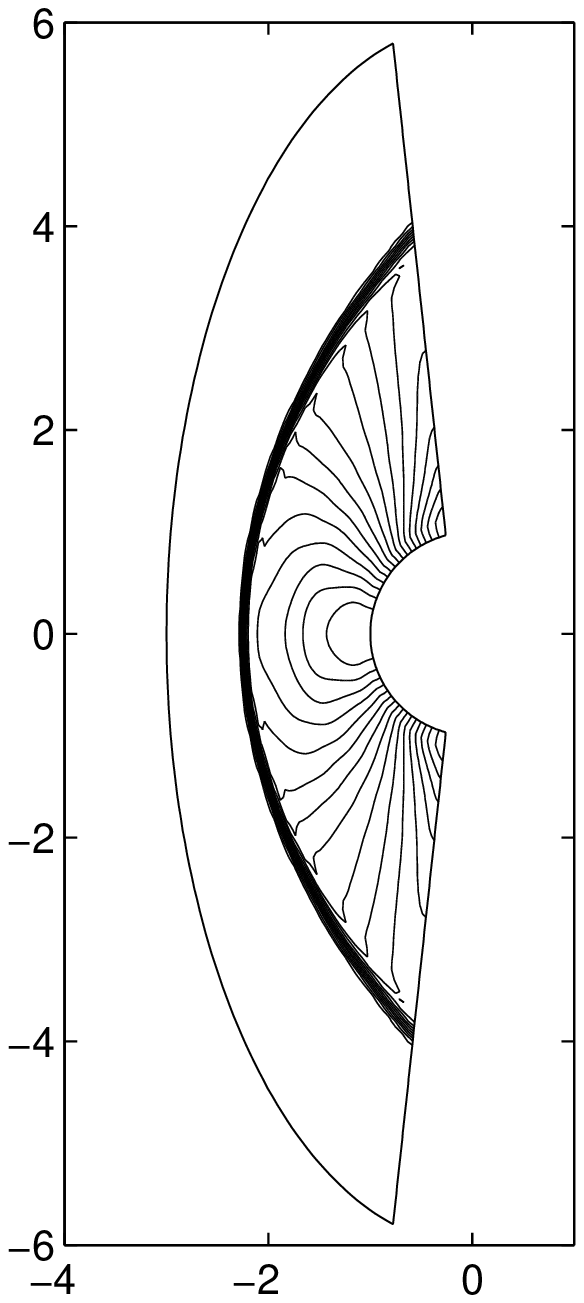}}
\caption{Flow pass a cylinder.
         (a) Physical grid; (b) Pressure with 20 contours.}
\label{cylin}
\end{figure}

The flow field is initialized with a Mach 3 shock moving from the inlet, i.e.,
$\xi=0$, toward the cylinder. The state in front of the
shock is given as $(\rho_r, u_r, v_r, p_r)=(1.4, 1.0, 0, 1.0)$. We impose
a reflective  boundary condition on the surface of the cylinder, i.e.,
$\xi=1$, and the outflow boundary condition at $\eta=0$ and $\eta=1$.
In our simulation, a uniform mesh of $65\times 129$ in the computation
domain is used. FIG. \ref{cylin}(a) depicts a diagram of the mesh
(drawing every other grid line) in the physical domain. Like the case on
a rectangular domain, we double the computational domain in both $\xi$ and
$\eta$ directions when we approximate the first order derivative of the
flux using the Fourier spectral method. The numerical result
given in terms of pressure contours (obtained at $t$=4.5) is shown in
FIG. \ref{cylin}(b). In summary, our spectral scheme can be easily applied
to the nonrectangular domain, as long as the domain can be smoothly
transformed to the rectangular one.

\section{Conclusion Remarks}

In this work, we introduce the use of discrete singular convolution (DSC) filters
in the Fourier pseudospectral method (FPM) for the simulation of hyperbolic
conservation law systems. The fast Fourier transform (FFT) is used as the basic
scheme, while DSC low-pass filters are used to suppress unphysical oscillations.
The DSC filters are implemented in either the
Fourier domain or the physical domain, depending on the physical problem
under consideration. The Fourier domain implementation is easy and cost efficient,
whereas the physical domain implementation allows more flexibility to handle
some special boundary conditions. A sensor technique developed in  the
conjugate filter oscillation reduction (CFOR) scheme \cite{guwei,weigu,ZhoWeijcp}
is adopted in the present work to appropriately activate the DSC filter.
The fourth-order Runge-Kutta scheme is utilized for the time integration.

Extensive numerical experiments are considered to validate the proposed
approach and to demonstrate its usefulness. Excellent numerical results
are obtained for all problems examined. The DSC filters  turn
out to be very effective in removing the Gibbs oscillations produced
by the spectral approximation while maintaining a high accuracy. Moreover,
it is also highlighted that these filters can dramatically improve the
stability of the spectral method for the long time integration.
The main features of proposed method are the follows:\\
(1) The basic scheme, the FPM is of spectral accuracy, which is important
    for problems that require not only capturing shock waves, but also resolving
    fine flow structures.
    \\
(2) The use of FFT endows the method high efficience, which is desirable for
    large scale problems in science and engineering.
    \\
(3) Numerical dissipation of the DSC filters is adjustable for a given problem,
    which endows the proposed method the ability to handle problems that are sensitive
    to numerical dissipation. The interaction of turbulence and shock is arguably a problem
    of this kind.
    \\
(4) The Fourier domain algorithm can be regarded as a windowed Fourier pseudospectral 
     method for shock-capturing.
    The alternative use of the DSC filter in the physical domain makes the proposed method
    feasible to problems requires special boundary conditions, e.g., reflective boundary
    condition.
\\
(5) The proposed scheme by-passes finding the characteristics which is required in many
    other shock-capturing schemes.
\\
(6) The use of the method in higher dimensions is straitforward.

It is note that the state of art shock-capturing schemes, such as WENO and
central schemes, do not depend on any adjustable parameter.
To make the proposed scheme robust and user friendly, future work will focus on
developing sophisticated sensors so that the filter parameter will be chosen
automatically for each given problem. We believe that this can be accomplished
by analyzing the Fourier power spectrum distribution of the approximate
solution at selected time steps.

{\bf Acknowledgments}
{This work was supported in part by Michigan State University.}

\newpage

\newpage

\newpage

\newpage

\begin{center}
\begin{table}
\begin{tabular}{ccccc}
        mesh               &  error  &   FFT     &   C4      &   DSC       \\ \hline
\multirow{2}{15mm}{$N$=32} &  $L_1$  &  5.58E-5  &  1.96E-3  &  1.27E-4  \\
                           &  $L_2$  &  1.27E-4  &  5.67E-3  &  4.07E-4  \\ \hline
\multirow{2}{15mm}{$N$=64} &  $L_1$  &  2.33E-8  &  1.32E-4  &  8.18E-8  \\
                           &  $L_2$  &  7.94E-8  &  4.39E-4  &  4.07E-7  \\ \hline
\multirow{2}{15mm}{$N$=128}&  $L_1$  &  4.01E-11 &  8.90E-6  &  3.99E-11  \\
                           &  $L_2$  &  5.09E-10 &  2.96E-5  &  5.09E-10  \\
\end{tabular}
\vspace{0.5cm}
\caption{Advective 2D isentropic vortex. $L_1$ and $L_2$ errors for the density
         at $t$=2. The CFL number is 0.01 for all schemes. C4 denotes the
         fourth-order accurate, conservative centered scheme.}
\label{vortextable1}
\end{table}
\end{center}

\vspace{2cm}

\begin{center}
\begin{table}
\begin{tabular}{cccccccc}
        mesh               &  error & $t$=100 & $t$=200 &  $t$=400& $t$=600 & $t$=800 & $t$=1000 \\ \hline
\multirow{2}{15mm}{$N$=64} &  $L_1$ & 4.63E-7 & 5.05E-7 & 1.00E-6 & 1.44E-6 & 2.07E-6 & 2.77E-6 \\
                           &  $L_2$ & 1.99E-6 & 1.23E-6 & 2.90E-6 & 4.42E-6 & 6.02E-6 & 7.59E-6 \\
\end{tabular}
\vspace{0.5cm}
\caption{Advective 2D isentropic vortex. $L_1$ and $L_2$ errors for the
         density at different times with $N$=64 (CFL=0.5, $\eta$=0.5).}
\end{table}
\end{center}

\newpage

\appendix
\begin{center}
  \textbf{APPENDIX}\\
  The optimal values ($r$) for DSC-RSK lowpass filters in our numerical experiments.
\end{center}

\begin{center}
\begin{table}
\begin{tabular}{ccc}
 No. of Example  &   case  &   $r$        \\   \hline
 \multirow{2}{5mm}{1} &  $N$=128  &  0.6  \\
                      &  $N$=256  &  0.8  \\   \hline
 \multirow{2}{5mm}{2} &  $N$=128  &  0.6 \\
                      &  $N$=256  &  0.8  \\   \hline
 \multirow{2}{5mm}{3} &  $N$=129  &  0.7  \\
                      &  $N$=257  &  0.8  \\   \hline
 \multirow{2}{5mm}{4} &  $N$=129  &  0.6  \\
                      &  $N$=257 &  0.6  \\   \hline
 \multirow{2}{5mm}{5} &  $N$=65   &  0.8  \\
                      &  $N$=129  &  0.8  \\   \hline
 \multirow{2}{5mm}{6} &  Lax      &  0.95 \\
                      &  Sod      &  1.1  \\   \hline
 \multirow{4}{5mm}{7} & $\kappa$=13  &  2.0  \\
                      & $\kappa$=26  &  2.0  \\
                      & $\kappa$=39  &  2.1  \\
                      & $\kappa$=52  &  2.1  \\ \hline
 \multirow{2}{5mm}{8} &  $N$=129  &  2.0 \\
                      &  $N$=257  &  2.1 \\   \hline
 \multirow{1}{5mm}{9} &           &  2.1  \\   \hline
 \multirow{1}{5mm}{10}&           &  2.8  \\   \hline
 \multirow{2}{5mm}{11}&  $\eta$=1.0  &  3.2  \\
                      &  $\eta$=0.5  &  2.8  \\ \hline
 \multirow{1}{5mm}{12}&           &  1.2
\end{tabular}
\end{table}
\end{center}

\end{document}